\documentclass[11pt]{amsart}
\usepackage{amsmath}
\usepackage{amsxtra}
\usepackage{amscd}
\usepackage{amsthm}
\usepackage{amsfonts}
\usepackage{amssymb}
\usepackage{eucal}
\usepackage[dvips]{color}
\newcommand{\nn}{\nonumber}
\newcommand{\bea}{\begin{eqnarray}}
\newcommand{\ena}{\end{eqnarray}}
\newcommand{\be}{\begin{eqnarray*}}
\newcommand{\en}{\end{eqnarray*}}
\def\bel{\begin{eqnarray}}
\def\enl{\end{eqnarray}}
\newcommand{\C}{{\mathbb C}}

\newcommand{\Z}{{\mathbb Z}}
\newcommand{\al}{{\alpha}}

\newcommand{\la}{{\lambda}}

\def\R{\mathcal R}
\numberwithin{equation}{section}

\numberwithin{equation}{section}
\newtheorem{thm}{Theorem}[section]
\newtheorem{prop}[thm]{Proposition}
\newtheorem{lem}[thm]{Lemma}
\newtheorem{cor}[thm]{Corollary}
\theoremstyle{remark}

\newtheorem{definition}[thm]{Definition}

\newtheorem{conj}[thm]{Conjecture}

\newcommand{\g}{\mathfrak{g}}
\newcommand{\h}{\mathfrak{h}}
\newcommand{\na}{\mathfrak{n}}

\newcommand{\hk}{\hookrightarrow}

\newcommand{\T}{\otimes}

\newcommand{\Uv}{U_v(\mathfrak{g})}
\newcommand{\Uvv}{U_{v^{-1}}(\mathfrak{g})}
\newcommand{\V}{\mathcal V}
\newcommand{\K}{{\mathbb K}}
\newcommand{\one}{{\bf 1}}
\newcommand{\e}{{\epsilon}}
\newcommand{\ga}{{\gamma}}
\newcommand{\id}{{\rm id}}
\newcommand{\End}{\mathop{\rm End}}
\newcommand{\bV}{\overline{\mathcal{V}}}

\newcommand{\cX}{\mathcal{X}}
\newcommand{\bcX}{\bar{\mathcal{X}}}
\newcommand{\W}{\mathcal{W}}
\newcommand{\bo}{B}
\renewcommand{\na}{N}
\newcommand{\Tor}{\rm Tor}
\begin{document}

\title[Fermionic formulas and difference Toda Hamiltonian]
{Fermionic formulas for eigenfunctions
of the difference Toda Hamiltonian}

\author{B. Feigin, E. Feigin, M. Jimbo, T. Miwa and E. Mukhin}
\address{BF: Landau Institute for Theoretical Physics,
Russia, Chernogolovka, 142432, prosp. Akademika Semenova, 1a,   \newline
Higher School of Economics, Russia, Moscow, 101000,  Myasnitskaya ul., 20 and
\newline
Independent University of Moscow, Russia, Moscow, 119002,
Bol'shoi Vlas'evski per., 11}
\email{bfeigin@gmail.com}
\address{EF:
Tamm Theory Division, Lebedev Physics Institute, Russia, Moscow, 119991,
Leninski pr., 53 and\newline
Mathematical Institute, University of Cologne, Weyertal 86-90, D-50931,
Cologne, Germany}
\email{evgfeig@gmail.com}
\address{MJ: Graduate School of Mathematical Sciences,
The University of Tokyo, Tokyo 153-8914;
Institute for the Physics and Mathematics of the Universe,
Kashiwa, Chiba 277-8582, Japan}
\email{jimbomic@ms.u-tokyo.ac.jp}
\address{TM: Department of Mathematics,
Graduate School of Science,
Kyoto University, Kyoto 606-8502,
Japan}\email{tetsuji@math.kyoto-u.ac.jp}
\address{EM: Department of Mathematics,
Indiana University-Purdue University-Indianapolis,
402 N.Blackford St., LD 270,
Indianapolis, IN 46202}\email{mukhin@math.iupui.edu}
%\date{\today}

\begin{abstract}
  We use the Whittaker vectors and the Drinfeld Casimir element to
  show that eigenfunctions of the difference Toda Hamiltonian can be
  expressed via fermionic formulas. Motivated by the combinatorics of
  the fermionic formulas we use the representation theory of the
  quantum groups to prove a number of identities for the coefficients
  of the eigenfunctions.
\end{abstract}

\maketitle

%%%%%%%%%%%%%%%%%%%%%%%%%%%%%%%%
\section{Introduction}\label{sec:1}

The goal of this paper is to derive fermionic formulas for
eigenfunctions of the finite difference Toda Hamiltonian $H_{Toda}$
and to study these fermionic formulas. Eigenfunctions of
$H_{Toda}$ have been studied recently in connection with quantum
cohomology of flag manifolds (see \cite{GiL}, \cite{BrFi}),
Whittaker vectors (see \cite{Br}, \cite{Sev}, \cite{Et}),
Macdonald polynomials and affine Demazure characters (see \cite{GLO}).
In particular, an important connection with the representation theory
of quantum groups was established.
In our paper we show how fermionic formulas naturally appear in the
representation-theoretical terms. On the other hand, these formulas
can be studied from purely combinatorial point of view.
In the paper we combine these two approaches.
We give some details below.

\subsection{Central elements and Whittaker vectors.}\label{CW}
As we have already mentioned, the representation theory of quantum groups
plays a very important role in the study of finite difference Toda Hamiltonian.
In particular, one can construct eigenfunctions of $H_{Toda}$ using
Whittaker vectors in Verma modules (see \cite{Br}, \cite{Sev}, \cite{Et}).
In this paper we use pairing
of Whittaker vectors with the dual ones.

Let $\g$ be a complex simple Lie algebra of rank $l$ and let
$U_v({\g})$  and $\Uvv$ be two quantum groups with parameters $v$ and $v^{-1}$.
Let $P$, $Q$  (resp. $P_+, Q_+$) be the weight and root lattices of $\g$
(resp. their positive parts) and
let
$\V^\la=\sum_{\beta\in Q_+} (\V^\la)_\beta$ and $\bV^\la=\sum_{\beta\in Q_+}(\bV^\la)_\beta$
be Verma modules of $U_v(\g)$ and $\Uvv$, respectively.
In order to define a Whittaker vector $\theta^\la$ in the completion
$\prod_{\beta\in Q_+} (\V^\la)_\beta$ of the Verma module
$\V^\la$ one fixes elements $\nu_i\in P$
and scalars $c_i$ ($1\le i\le l$).
Then the Whittaker vector, associated with these data, is defined by the condition
\bea\label{whit1}
E_i K_{\nu_i}\ \theta^\la
=\frac{c_i}{1-v^2}\ \theta^\la
\ena
(for simplicity, in Introduction, we assume that $\g$ is simply-laced).
Here $E_i\in U_v({\frak g})$ are the Chevalley generators (which act as annihilating operators) and
$K_{\nu_i}$ are certain elements from the Cartan subalgebra, associated with $\nu_i$.
Similarly, one defines the dual Whittaker vector $\bar\theta^\la$ in
the completion  of $\bV^\la$ by the formula
\bea\label{whit2}
\bar E_i \bar K_{\nu_i}\ \bar \theta^\la
=\frac{c_i^{-1}}{1-v^{-2}}\ \bar \theta^\la
\ena
The central object of our paper is the following function
\be
J^\la_\beta=v^{-(\beta,\beta)/2+(\la,\beta)}\
(\theta^\la_\beta,\bar\theta^\la_\beta),
\en
where $\theta^\la_\beta\in (\V^\la)_\beta$ is the weight $\la-\beta$ component of the Whittaker
vector and $(~,~)$  is the natural non-degenerate pairing between $\V^\la$ and $\bV^\la$.
It can be shown that $J^\la_\beta$ is independent of possible choices of $\nu_i$ and $c_i$.

Consider the generating function
\be
F(q,z_1,\dots,z_l,y_1,\dots,y_l)=\sum_\beta J^\la_\beta \prod_{i=1}^l y_i^{(\beta,\omega_i)},
\en
where $z_i=q^{-(\la,\al_i)}$, $q=v^2$ and   $\omega_i$ (resp. $\alpha_i)$
are fundamental weights (resp. simple roots).
Then $F$ is known to be an eigenfunction of the quantum difference Toda
operator (\cite{Sev}, \cite{Et}). In order to prove this statement
one uses central elements of the quantum group.
Roughly, the procedure works as follows.
If $u$ is a central element, then the scalar product
\bea\label{u}
(u\theta^\la_\beta,\bar\theta^\la_\beta)
\ena
can be written in two ways. On the one hand, one can compute the action of $u$
on $\V^\la$ (the corresponding  scalar). On the other hand, if a precise
formula for $u$ is known then one can compute $(\ref{u})$ using the relation
\be
(F_iw,\bar w)=(w,\bar E_i\bar w)
\en
and formulas $(\ref{whit1})$, $(\ref{whit2})$.

The Toda Hamiltonian appears when one uses the central element written as
the trace of products of $R$ matrices in finite-dimensional $U_v(\g)$ modules.
Our key observation is that if
the Drinfeld Casimir element is used instead then one obtains
a recursion relation for $F$ which leads to the fermionic formulas.
In the next subsection we describe those formulas in more details.

\subsection{Fermionic formulas.}
Fermionic formulas appear in different problems of representation theory and
mathematical physics (see for example  \cite{BM}, \cite{FJMMT}, \cite{HKOTT}, \cite{SS}).
Let us describe the class of formulas we treat in our paper.

Let $[r,s]=\{t\in\Z\mid r\leq t\leq s\}$ be a subset of $\Z$, where $r,s$ are integers
or $\pm\infty$. Let $V$ be a vector space with a basis $e_{i,t}$ labeled by
pairs $1\le i\le l$, $t\in [r,s]$. Let
$\Gamma_+=\{\sum_{(i,t)}m_{i,t} e_{i,t}| m_{i,t}\in\Z_{\geq0}\}$ be the
positive part of the lattice generated by $\{e_{i,t}\}$.
We fix a quadratic form $\langle\cdot,\cdot\rangle$ on $V$ and a vector $\mu\in V$.
Further, define  maps $w$ and $d$ from $V$ to the $l$-dimensional vector space
with a basis $p_1,\dots, p_l$ via the formulas
\be
w(\sum_{(i,t)} m_{i,t}e_{i,t})=\sum_{i=1}^l p_i \sum_{t\in [r,s]} m_{i,t},  &&
d(\sum_{(i,t)} m_{i,t}e_{i,t})=\sum_{i=1}^l p_i \sum_{t\in [r,s]} tm_{i,t}.
\en
Define functions $I_m$ depending on $q$, $z=(z_1,\dots,z_l)$ and
$m=(m_1,\dots,m_l)$ as follows
\bea\label{I_m}
I_m(q,z)=\sum_{w(\gamma)=m}
z^{d(\gamma)}
\frac{q^{\langle \gamma,\gamma\rangle + \langle\mu,\gamma\rangle}}{(q)_\gamma},
\ena
where the summands are labeled by $\gamma=\sum_{(i,t)} m_{i,t}e_{i,t}\in\Gamma_+$ and
$(q)_\gamma=\prod_{(i,t)} (q)_{m_{i,t}}$, $z^{d(\gamma)}=\prod_{i=1}^l z_i^{d(\gamma)_i}$.
We call the right hand side of $(\ref{I_m})$ a fermionic formula.
The generating function $F(q,z,y)=F(q,z_1,\dots,z_l,y_1,\dots,y_l)$ is given by the formula
\bea\label{genfun}
F(q,z,y)=\sum_{m}\ y^m I_m(q,z),\quad y^m=y_1^{m_1}\dots y_l^{m_l}.
\ena

Let the matrix of the quadratic form
$\langle\cdot,\cdot\rangle$ be a tensor product
$D=C\otimes G(r,s)$, where $C$ is the Cartan matrix of $\g$ (we assume
here that $C$ is symmetric) and
$G=(G_{t,t'})_{i,j\in[r,s]}$, $G_{t,t'}=\min(t,t')$.
Such matrices appear
in \cite{DS}, \cite{S} in the fermionic
formulas for the Kostka polynomials.
Let $[r,s]=[0,\infty).$
Then functions $I_m(q,z)$ satisfy the following recursion relation:
\bea\label{Ima}
I_m(q,z)=\sum_{0\le a\le m}
\frac{z^aq^{W(a)}}{(q)_{m-a}}I_a(q,z),
\ena
where $W(a)=\frac{1}{2}(Ca\cdot a-\mathrm{diag}C\cdot a)$, $\cdot$ denotes the
standard scalar product and $0\le a\le m$ abbreviates the set of
inequalities $0\le a_i\le m_i$.
The relation $(\ref{Ima})$ shows that $I_m(q,z)$ are
determined by $I_0(q,z)$.

Recall the functions $J^\la_\beta$. Using the Drinfeld Casimir element and the procedure described
in the end of subsection \ref{CW}, we show that $J^\la_\beta$
satisfy the relation
\be
J^\la_\beta=\sum_{\beta'}
\frac{1}{(q)_{\beta-\beta'}}
q^{(\beta',\beta')/2-(\la+\rho,\beta')}J^{\la}_{\beta'}.
\en
This leads to the following identification
\be
J^\la_{\beta}=%I_m(q),
I_m(q,z),
\ \beta=\sum_i {m_i}\al_i,\ z=q^{-(\la,\al_i)}.
\en
In particular, this gives a fermionic formula for eigenfunctions of
$H_{Toda}$.

Fermionic sums can be considered as a sort
of statistical sum for some ``models".
The models depend on parameters $r$ and $s$ and enjoy many ``physical''
combinatorial properties. For example, we
look into what %is happening
happens with the fermionic
sums when the parameters,
e.g., $r$ and $s$, go to infinity. We sort the terms by the dependence
on the parameters which go to infinity and call the result the quasi-classical decomposition.
Then we expect that the quasi-classical decompositions
are exact, which means that the
coefficients in a decomposition are summed up to rational functions
and the result gives correct
formulas for finite values of parameters.
That expectation predicts recursion
relations for the fermionic sums $I_m$.

We then prove the recursion relations using the Whittaker vectors and
the representation theory of the  quantum group $U_v(\g)$,  (see Theorems
\ref{thm:XJ},  \ref{thm:k=0-1},  \ref{thm:J[0,k]}).
In some cases the relations become finite.
From the point of view of fermionic sums it means
the vanishing property: some fermionic expressions
must be zero.  The quantum group approach explains
the vanishing property as well.
Namely, some terms in the recursions are
zero because the corresponding weight space is zero in the
irreducible representation of the quantum group.

It is well known  (see \cite{Et})
that the eigenfunctions of the difference Toda
Hamiltonian can be obtained as a certain limit of the Macdonald
polynomials. In the case of $\frak{sl}_n$, one of the recursions we
prove (see $(\ref{id1})$) is the corresponding limit of the Pierri rule
for the Macdonald polynomials. It is interesting to study the other
identities in relation with the Macdonald polynomials. We hope to
address this problem in future publications.

\subsection{Affine Lie algebras: motivations and further directions.}
In this subsection we discuss various connections between eigenfunctions
of $H_{Toda}$ and representation theory of affine Kac-Moody algebras. The fermionic formulas
provide a very useful tool for computation of various "affine" characters.
Though we do not treat this subject in the main body of the paper,
it was our original motivation for studying
the  "fermionic part" of the quantum difference Toda story.
We are not providing any proofs here. We hope to
return to this subject in more details elsewhere.

\subsubsection{Refined characters.}
Let $\frak a$ be a Lie algebra and $W$ its representation.
The character of $W$ is an expression
$$
\chi(z_1,\ldots,z_n)
={\rm tr}_W\left(z_1^{a_1} \cdots z_n^{a_n}\right),
$$
where $a_i$ are some commuting elements of $\frak a$.
Usually, $\frak a$ is semi-simple and
$(a_1,\ldots,a_n)$ is a basis of
its Cartan subalgebra. There is a simple way to ``refine" the
character.
To do it,
suppose that
$W$ is a cyclic representation with a cyclic vector $v$,
and choose a subspace
$S\subset U({\frak a})$ such that $1\in S$.
We define subspaces $F_j\hk W$, $j=0,1,\dots$:
$F_0=\C\cdot v$ and
$F_{j+1}=S\cdot F_j$.
Assuming that $F_j$ converge to $W$ we obtain a filtration in $W$.
Suppose also that the space $S$
is invariant with respect to the adjoint action of $a_i$ for all $i$: $$[a_i,S]\subset S.$$
In this case, if $F_0$ is $\{a_i\}$ invariant, i.e., $a_iF_0\subset F_0$, then all spaces $F_j$
in the filtration are $\{a_i\}$-invariant.
Consider now the associated graded space
$$
\overline W=F_0\oplus(\bigoplus_{j>0} F_j/F_{j-1}).
$$
On $\overline W$ we have an action of
$\{a_i\}$ and also the action of an additional operator $b$,
which acts as  the constant $j$ on $F_j/F_{j-1}$.
Now, we define the refined character
\begin{align*}
\chi(z_1,\ldots,z_n;y)&
={\rm tr}_{\overline W}
\left(z_1^{a_1}\cdots z_n^{a_n} y^b\right),\\
\chi(z_1,\ldots,z_n;1)&=\chi(z_1,\ldots,z_n).
\end{align*}

We consider the case
$\frak a=\frak h\oplus \widehat{\frak n}$.
Here $\widehat{\frak n}={\frak n}\otimes\C[t,t^{-1}]$,
$\frak n$ is a maximal nilpotent subalgebra of
a finite-dimensional semi-simple Lie algebra $\frak g$,
$\tilde{\mathfrak h}\simeq \frak{h}\T 1\oplus\C d$  is the Cartan subalgebra of
$\frak g\otimes\C[t,t^{-1}]\oplus \C d$ and $d=xd/dx$ is the grading operator.
Let $e_1,e_2,\ldots,e_l$ be the generators of $\frak n$ and
$e_i[j]=e_i\T t^j$
be the corresponding generators of $\widehat{\frak n}$.
We define currents $e_i(x)=\sum_{j\le 0} e_i[j]x^{-j}$.
Fix a basis
%$\{L_0,h_1,\ldots,h_l\}$
$\{d,h_1,\ldots,h_l\}$
in $\tilde{\mathfrak h}$.
As a representation $W$ we take the induced module
generated by the vacuum vector $v$ satisfying
$e_i[j]v=0$ for $j>0$.
We have the character
$$
\chi(q,z_1,\ldots,z_l)
=
%{\rm tr}_W\bigl(q^{L_0}z_1^{h_1}\cdots z_l^{h_n}\bigr).
{\rm tr}_W\bigl(q^{-d}z_1^{h_1}\cdots z_l^{h_l}\bigr).
$$

\subsubsection{The $A_1$ case.}
We consider  $\g=\frak{sl}_2$. In this case $\frak{n}$ is
one-dimensional and is spanned by an element $e$. We fix $W$ to
be an induced $\frak{a}$ module with a cyclic vector $v$ satisfying
$\frak{h}\cdot v=0$ and $e[j]v=0$, $j>0$. Let $S$ be the subspace
spanned by coefficients of the expansion of $e(x)^s$, $s\geq0$
as series in $x$. Then we obtain
\begin{align*}
\chi(q,z)&=\frac1{(z)_\infty},\\
\chi(q,z;y)&=\sum_{m\geq0}\frac{y^mz^mq^{m^2}}{(q)_m(z)_m},
\end{align*}
where $(z)_m=\prod_{n=1}^m(1-q^{n-1}z)$.
We prove in Appendix A (see \eqref{f r in F G}) that
$\chi(q,z;y)$ differs from the generating function of $J^\la_\beta$
(or, equivalently, of $I_m$) by a simple factor.
We now give fermionic expression for the quantity
$\chi(q,z;y)$.

Introduce an algebra generated by Fourier coefficients
$c_j[s]$ of the currents
$c_j(x)=\sum_{s\le 0}c_j[s]x^{-s}$,
$j=0,1,2,\ldots$. The defining relations are
$[c_{j_1}[s_1],c_{j_2}[s_2]]=0$
for all $j_1,j_2,s_1,s_2$, and $c_j(x)^2=0$.
Note that $c_j(x)$ can be constructed as vertex operators with
momentum $p_j$ such that $\langle p_j,p_j\rangle=2$
and $\langle p_i,p_j\rangle=0$ for $i\not=j$.
Let
$$
e(\varepsilon,x)=\sum_j\varepsilon^jc_j(x),
$$
where $\varepsilon$ is a formal variable.
Let $\mathcal A$ be an algebra over
$\C[[\varepsilon]]$, the ring of formal
power series in $\varepsilon$, which
is generated by $e(\varepsilon,x)$.
In $\mathcal A$ there exists a subspace
$S$ spanned by the elements $\{e(\varepsilon,x)^n;n\geq0\}$,
and a filtration $\mathcal F_j$
such that $\mathcal F_0=\C\cdot1$ and $\mathcal F_j=\mathcal F_{j-1}+ S\cdot\mathcal F_{j-1}$.
The associated graded space
$\overline{\mathcal A}$ naturally has a
structure of commutative algebra.
It is generated by the space
$\mathcal F_1/\mathcal F_0$.
Actually, $\overline{\mathcal A}$ is a quadratic algebra.
It is a free module over
$\C[[\varepsilon]]$.
Let us consider the specialization $\overline{\mathcal A}_0$ at $\varepsilon=0$.
The algebra $\overline{\mathcal A}_0$ is
generated by the currents $d_1(x)=c_0(x)$,
$d_2(x)=c_0(x)c_1(x)$, $d_3(x)=c_0(x)c_1(x)c_2(x)$, etc..
The defining relations in
$\overline{\mathcal A}_0$ are quadratic
and takes the form involving
derivatives of the currents:
\begin{align}
%d_i(x)\cdot d^{(l)}_j(x)=0,\quad 0\leq l\leq 2\min(i,j).
d_i(x)\cdot d_j(x)^{(l)}=0,\quad 0\leq l\leq 2\min(i,j).
\label{defining relation}
\end{align}

The representation $W$ for the algebra
$\mathcal A=\{e(\varepsilon,x)\}$ is also defined
on the
free module
over the ring $\C[[\varepsilon]]$,
and after substituting $\varepsilon=0$ we get
a representation over $\overline{\mathcal A}_0$.
It is quadratic with simple relations, and
this construction gives us a fermionic formula for the
refined character $\chi$:
$$
\chi(q,z;y)=\sum_{\{m_j\}}
\frac{q^{2\sum\min(i,j)m_im_j}z^{\sum jm_j}y^{\sum m_j}}
{\prod(q)_{m_j}}.
$$

\subsubsection{General $\g$ and multi-filtration.}
In order to treat the general case we need to replace the filtration
$F_i$ by a ``multi-filtration".
Let $S_i$, $i=1,\ldots,l$,
be the subspaces of $U(\widehat{\frak n})$
spanned by coefficients of the expansion of $e_i(x)^s$, $s\geq0$,
as series in $x$,
and let $R_i=\oplus_{k\leq i}S_k$, where $R_0=S_0=\C\cdot 1$.
Note that
$R_1\subset R_2\subset\cdots\subset R_l$.
Then, we define
\begin{align*}
&F_{j_1}=R_1\cdot F_{j_1-1},\quad F_0=\C\cdot v,\\
&F_{j_1,j_2}=R_2\cdot F_{j_1,j_2-1},\quad
F_{j_1,0}=F_{j_1},\\
&\cdots.
\end{align*}
We define
\begin{align*}
&\overline W=\oplus_{j_1,\ldots,j_l}\overline
W_{j_1,\ldots,j_l},\\
&\overline W_{j_1,\ldots,j_l}=
F_{j_1,\ldots,j_l}/\sum_iF_{j_1,\ldots,j_i-1,\ldots,j_l}
\end{align*}
and denote by $b_i$ the operator which
gives the constant $j_i$ on
$\overline W_{j_1,\ldots,j_l}$.

In this way, we get  the subspaces
%$W_{j_1,\ldots,j_l}$
$F_{j_1,\ldots,j_l}$
and the associated graded
spaces
%$\overline W$;
${\overline W}_{j_1,\ldots,j_l}$;
the latter have the actions of the grading
operators $b_1,\ldots,b_l$.
Since the operators %$L_0,h_1,\ldots,h_l$
$d,h_1,\ldots,h_l$
act on $\overline W$ in an evident way,
we can write
$$
\chi(q,z_1,\ldots,z_l;y_1,\ldots,y_l)
={\rm tr}_{\overline W}
%\left(q^{L_0}z_1^{h_1}\cdots z_l^{h_l}
\left(q^{-d}z_1^{h_1}\cdots z_l^{h_l}
y_1^{b_1}\cdots y_l^{b_l}\right).
$$
We conjecture that this refined character gives an eigenfunction of
the conjugated quantum Toda Hamiltonian (see \cite{BrFi}, \cite{GiL}, Appendix A).
We now explain how to obtain a fermionic formula in the general case.

\subsubsection{General fermionic formula.}
Let $L_1$ be the vacuum representation of $\widehat{\frak g}$
of level $1$ with highest weight vector $v$,
and let $V$ be the principal subspace in $L_1$, i.e.,
$V=U(\widehat{\frak n})\cdot v$.
The space $V$, as a representation of $U(\widehat{\frak n})$
can be described by using the space $W$ as
$$V=W/\sum_ie_i(x)^2W.$$
Consider the tensor product of infinitely many copies of $V$ labeled
by $j=0,1,2\ldots$, and denote by
$a^{(j)}_i$ the operator $e_i$ acting on the $j$-th copy:
$$
a^{(j)}_i=1\otimes\cdots
\otimes\underbrace{e_i}_{\textstyle j\hbox{-th}}\otimes \cdots.
$$
Now, set
$$
e_i(\varepsilon,x)=\sum_{\alpha\geq0}\varepsilon^\alpha a^{(\alpha)}_i(x).
$$
After substituting $\varepsilon=0$ we get an algebra generated
by
$c_i^{(\beta)}(x)=a^{(0)}_i(x)a^{(1)}_i(x)\cdots
a^{(\beta-1)}_i(x)$.
The currents $c_i^{(\beta)}(x)$
generate an algebra with the quadratic relations
\begin{align}
c_i^{(\beta)}(x)
c_j^{(\gamma)}(x)^{(l)}=0
\hbox{ if }l\leq\min(\alpha,\beta)\cdot C_{i,j}
\end{align}
where $(C_{i,j})$ is the Cartan matrix of
$\frak g$. As a result we have a fermionic formula
for the refined character of $W$:
\bea\label{F formula}
&&\chi(q,z_1,\ldots,z_l;y_1,\ldots,y_l)
=
\sum_{\{ n^{(t)}_i\}}
\frac{q^{B}}{\prod_{t\ge0}\prod_{i=1}^l(q_i;q_i)_{n^{(t)}_i}}
\prod_{i=1}^lz_i^{\sum_{t\ge0}tn^{(t)}_i}
\prod_{i=1}^ly_i^{\sum_{t\ge0} n^{(t)}_i},
\\
&&
B=\sum_{t\ge0}t
\left(\sum_{i,j=1}^l\frac{1}{2}b_{i,j}n^{(t)}_in^{(t)}_j
-\sum_{i=1}^l d_in^{(t)}_i\right)
+\sum_{t<t'} t\sum_{i,j=1}^l b_{i,j}n^{(t)}_in^{(t')}_j,
\nn
\ena
where $b_{i,j}=(\al_i,\al_j)$, $d_i=b_{i,i}/2$
and $q_i=q^{d_i}$. We note that if $\g$ is simply laced the formula
\eqref{F formula} corresponds to \eqref{I_m} and \eqref{genfun}.

\subsection{Plan of the paper}
Now let us outline the content of our paper.\\
In Section 2 we introduce the
fermionic sums. We study quasi-classical limits of such
formulas and various recursion relations.\\
In Section 3 we prove the fermionic
formulas for the scalar products of the
Whittaker vectors with dual ones. We also discuss the general procedure
(based on the center of the quantum group),
which produces equations satisfied by $J^\la_\beta$.\\
In Section 4 we study the quasi-classical
decompositions using the representation theory of $U_v(\frak g)$.
We prove the recursion relations and vanishing properties from Section 2.\\
In $A_l$ case there is an alternative simple
way to prove that the fermionic formula satisfies
the Toda equation. We give this proof in Appendix A.\\
In Appendix B,
we prove a proposition on the singular vectors in the tensor
product of two Verma modules.
We need this lemma to prove the vanishing property.

%%%%%%%%%%%%%%%%%%%%%%%%%%%%%%%%%%

\section{Fermionic sums}\label{sec:2}
\subsection{Fermionic sums on a finite interval}\label{finite interval}
Let $l\in\Z_{\geq1}$ be a positive integer, and let $[r,s]=\{t\in\Z\mid r\leq t\leq s\}$
be a finite interval in $\Z$. Here $r,s$ are integers, but in later subsections
we also consider the case $r=-\infty$ and/or $s=\infty$.

Let $C=(C_{i,i'})_{1\leq i,i'\leq l}$ be a symmetric matrix,
and let $\mu=(\mu_{i,t})_{1\leq i\leq l,t\in[r,s]}$ be a vector.
Let $m=(m_1,\ldots,m_l)$ be a set of non-negative integers.
We will define a fermionic sum $I_{C,\mu,m}(q,z|r,s)$ on the interval $[r,s]$
corresponding to the data $C,\mu$ and $m$.
Each fermionic sum is a rational function in $q$ and $z=(z_1,\ldots,z_l)$. It is defined as a sum of rational functions
parameterized by a configuration of particles ${\bf m}$, which we will explain below.

We call a tuple of non-negative integers \be {\bf m}=\{m_{i,t}:\ 1\le
i\le l, r\le t\le s\} \en a configuration of particles. An integer
$i\in[1,l]$ is called color and $t\in[r,s]$ is called weight of the
particle.  The non-negative integer $m_{i,t}$ represents the number of
particles with color $i$ and weight $t$. For a configuration of
particles ${\bf m}$, we associate a vector $\overline{\bf
  m}=(m_i)_{1\leq i\leq l}$ by \bea\label{m}
m_i=\sum_{t\in[r,s]}m_{i,t}.  \ena The number $m_i$ is the number of
particles with color $i$.

We set $(w)_n=\prod_{i=1}^n(1-q^{i-1}w)$ and  define
\be
(q)_{\bf m}=\prod_{(i,t)\in [1,l]\times [r,s]} (q)_{m_{i,t}}
\en
for any tuple ${\bf m}=(m_{i,t})$. We also use the standard scalar product
$({\bf m},{\bf n})=\sum_{(i,t)\in [1,l]\times [r,s]} m_{i,t}n_{i,t}$.

\begin{definition}
Let $m=(m_i)\in Z_{\ge 0}^l$ and $\mu=(\mu_{i,t})\in\Z^{l(s-r+1)}$ be two vectors.
A fermionic sum $I_{C,\mu,m}(q,z|r,s)$ is a function in $q$ and $z=(z_1,\dots,z_l)$
defined by
\bea
I_{C,\mu,m}(q,z|r,s)&=\sum_{\overline{\bf m}=m}\prod_{i=1}^l
z_i^{\sum_{t=r}^stm_{i,t}}
\frac{q^{Q_C({\bf m})+(\mu,{\bf m})}}{(q)_{\bf m}},\label{J}
\ena
where
\bea
Q_C({\bf m})&=\frac{1}{2}\left\{((C\otimes G){\bf m},{\bf m})-
({\rm diag}(C\otimes G), {\bf m})\right\}\label{Q},
\ena
the matrix $G$ is defined by
\be
G=(G_{t,t'})_{t,t'\in[r,s]},\qquad  G_{t,t'}=\min(t,t'),
\en
and for a matrix $X=(X_{i,j})$,  ${\rm diag}(X)$ signifies the vector consisting of
diagonal entries $X_{i,i}$.
\end{definition}
The quantity $I_{C,\mu,m}(q,z|r,s)$ is a Laurent polynomial in $z=(z_1,\ldots,z_l)$ with
coefficients which are Laurent polynomials in $q^{C_{i,i'}},q^{\mu_{i,t}}$ and rational functions in $q$.

We define the following formal power series in $y=(y_1,\ldots,y_l)$.
\bea
F_{C,\mu}(q,z,y|r,s)&=\sum_{m\in \Z^l_{\ge 0}}y^mI_{C,\mu,m}(q,z|r,s).\label{F}
\ena
Here we use the notation $y^m=\prod_{i=1}^ly_i^{m_i}$ for $l$-component vectors $y$ and $m$.
We use the convention
$F_{C,\mu}(q,z,y|r,s)=1$ if $r=s+1$.
We denote the functions $F_{C,\mu}(q,z,y|r,s)$, $I_{C,\mu,m}(q,z|r,s)$ in the case of $\mu=0$ by
$F_C(q,z,y|r,s)$, $I_{C,m}(q,z|r,s)$,
dropping the parameter $\mu$.

In what follows we also need another parametrization of configurations of
particles. Namely, with each ${\bf m}=(m_{i,t})$, $1\le i\le l$, $t\in[r,s]$,
we associate
the vector
\be
{\bf p}=(p_{i,j}),\quad 1\le i\le l, 1\le j\le m_i
\en
defined by
two conditions:
\begin{itemize}
\item $p_{i,1}\le p_{i,2}\le\cdots\le p_{i,m_i},\ 1\le i\le l$,
\item $m_{i,t}=\#\{j:\ p_{i,j}=t\}.$
\end{itemize}
It is easy to see that the correspondence ${\bf m}\leftrightarrow {\bf p}$ is one-to-one.
In the following lemma we rewrite powers of $z$ and $q$ in
\eqref{J} in terms of ${\bf p}$. To avoid confusion we
denote the function $Q_C(\bf m)$ written in $\bf p$ coordinates by
$\overline Q_C(\bf p)$.
\begin{lem}
We have
\bea
z_i^{\sum_{t=r}^stm_{i,t}}&=&z_i^{\sum_{j=1}^{m_i}p_{i,j}},\label{power z}\\
Q_C({\bf m})=\overline Q_C(\bf p)&=&
\frac{1}{2}\left(\sum_{(i,j),(i',j')}C_{i,i'}\min(p_{i,j},p_{i',j'})
-\sum_{(i,j)} C_{i,i} p_{i,j}\right).\label{quadratic}
\ena
\end{lem}

If we shift the parameters $\mu_{i,t}$ to $\mu_{i,t}+t\nu_i+\kappa_i$,
the sum $I_{C,\mu,m}(q,z|r,s)$ responds by a simple prefactor and $q$-shifts of $z_i$:
\begin{align}
I_{C,\mu+t\nu+\kappa,m}(q,z|r,s)=q^{\kappa\cdot m}I_{C,\mu,m}(q,q^\nu z|r,s).\label{shift}
\end{align}
Here, $\nu=(\nu_i)_{i\in[1,l]}$ and $(\mu+t\nu+\kappa)_{i,t}=\mu_{i,t}+t\nu_i+\kappa_i$.
We use the abbreviated notation $\kappa\cdot m=\sum_{i=1}^l\kappa_im_i$ and $q^\nu z=(q^{\nu_i}z_i)_{i\in[1,l]}$.

If $\mu=0$ and we shift the interval, we get a simple factor:
$$
I_{C,m}(q,z|r+k,s+k)=\left(z^mq^{W_{C,m}}\right)^kI_{C,m}(q,z|r,s),
$$
where
$$W_{C,m}=\frac12(Cm\cdot m-{\rm diag}\,C\cdot m).$$

\subsection{Fermionic sums on a semi-infinite interval}
We replace the interval $[r,s]$ in the construction in Section \ref{finite interval}
by an infinite interval. Let us consider the case $[r,s]=[0,\infty)$.
We use the abbreviation
\bea
I_{C,\mu,m}(q,z)=I_{C,\mu,m}(q,z|0,\infty),
\quad
I_{C,m}(q,z)=I_{C,0,m}(q,z|0,\infty).
\label{ICm}
\ena
In this case, the sum in the right hand side of \eqref{J} becomes an infinite sum.
We are interested in the case where the parameters
$\mu=(\mu_{i,t})_{(i,t)\in [1,l]\times [0,\infty)}$ are specialized so that the infinite sum is well-defined
as a rational function in $q,z_i,q^{C_{i,i'}},q^{\mu_{i,t}}$.

\begin{lem}\label{rational}
Suppose that there exists $t_0\geq0$ and $\nu_i,\kappa_i$ $(1\leq i\leq l)$ such that
\begin{eqnarray}
\mu_{i,t}=t\nu_i+\kappa_i\hbox{ if }t\geq t_0.\label{nu}
\end{eqnarray}
Then, the sum \eqref{J} is a polynomial in $q^{\mu_{i,t}}$ $((i,t)\in [1,l]\times [0,t_0-1])$
and $q^{\kappa_i}$ $(1\leq i\leq l)$
with coefficients which are rational functions in $q,z_i,q^{C_{i,i'}},q^{\nu_i}$ $(1\leq i,i'\leq l)$.
\end{lem}

The proof is easy, and we omit it.

There are simple relations between the fermionic sums over semi-infinite intervals with $\mu=0$.
Namely, we can reduce all cases to $[0,\infty)$. We have
\begin{align}
I_{C,m}(q,z|k,\infty)&=(z^mq^{W_{C,m}})^kI_{C,m}(q,z),\label{J1}\\
I_{C,m}(q,z|-\infty,k)&=(z^mq^{W_{C,m}})^kI_{C,m}(q,z^{-1}q^{-Cm+{\rm diag}\,C}).\label{J2}
\end{align}
Here we used the abbreviation $z^{-1}q^m=(z_1^{-1}q^{m_1},\ldots,z_l^{-1}q^{m_l})$.

The following proposition determines the fermionic sums recursively.
\begin{prop}\label{sec2:frec}
The rational functions $I_{C,m}(q,z)$ satisfy the recursion
\begin{align}
I_{C,m}(q,z)=\sum_{0\le a\le m}
\frac{z^aq^{W_{C,a}}}{(q)_{m-a}}I_{C,a}(q,z).\label{fermionic recursion}
\end{align}
The solution of this recursion is unique if we fix $I_{C,0}(q,z)=1$.
\end{prop}
\begin{proof}
We subdivide the fermionic sum in the left hand side into parts labeled by
$a=(a_1,\ldots,a_l)\in\Z_{\ge 0}^l$,
with $a_i$ being the number of color $i$
particles whose weights are larger than or equal to $1$.
Then this corresponds to the sum in the right hand side
because of \eqref{J1}.
The uniqueness is clear because the equation can be written as
$$(1-z^mq^{W_{C,m}})I_{C,m}(q,z)=\sum_{0\leq a<m}\frac{z^aq^{W_{C,a}}}{(q)_{m-a}}I_{C,a}(q,z).$$
\end{proof}
We call \eqref{fermionic recursion} the fermionic recursion.

In what follows we study functions $I_{C,\mu,m}(q,z)$ for  some special
values of $\mu$ of the type \eqref{nu}. We describe
$\mu$ in terms of corners and angles.

We say $\mu$ has a corner at $t\in[1,\infty)$
if the vector $\mu[t]$ given by $\mu_i[t]=\mu_{i,t+1}+\mu_{i,t-1}-2\mu_{i,t}$
is not zero. We call $\mu[t]$ the angle of $\mu$ at $t$.
We define $\mu_i[0]=\mu_{i,1}-\mu_{i,0}$,
and call $\mu[0]$ the angle of $\mu$ at $0$.
We say $\mu$ has a corner at $0$ if $\mu[0]\not=0$.
Let us discuss simple cases.

First, consider the case where
$\mu_{i,t}=t\nu_i+\kappa_i$ for all $t\geq0$. It reduces to the basic case
$I_{C,m}(q,z)$ by \eqref{shift}.

Second, we define the case with one corner. Namely, we consider
\begin{align}
\mu_{i,t}=
\begin{cases}
0&\hbox{ if }t\leq k;\\
(t-k)\nu_i&\hbox{ if }t\geq k.
\end{cases}\label{corner}
\end{align}
We denote the fermionic sum corresponding to this $\mu$ by $J^{(k,\nu)}_{C,m}(q,z)$.

\begin{lem}\label{J=II}
There exist mutual relations between $I_{C,m}(q,z),J^{(k,\nu)}_{C,m}(q,z),I_{C,m}(q,z|0,k)$.
We have
\begin{align}
J^{(0,\nu)}_{C,m}(q,z)&=I_{C,m}(q,q^\nu z),\label{relation 1}\\
J^{(k,\nu)}_{C,m}(q,z)&=
\sum_{0\le a\le m}\left(z^aq^{W_{C,a}}\right)^kI_{C,m-a}(q,q^{Ca}z|0,k-1)I_{C,a}(q,q^\nu z),
\label{relation 2}\\
&=\sum_{0\le a\le m}q^{a\cdot\nu}\left(z^aq^{W_{C,a}}\right)^{k+1}
I_{C,m-a}(q,q^{Ca}z|0,k)I_{C,a}(q,q^\nu z).\label{relation 3}
\end{align}
\end{lem}
\begin{proof}
%The proof is straightforward.
The relation \eqref{relation 1} is trivial.
In order to prove the other two relations we
cut the interval  of weights $[0,\infty)$ into two parts.
For \eqref{relation 2} these parts are $[0,k-1]$ and $[k,\infty)$ and
for \eqref{relation 3} they are $[0,k]$ and $[k+1,\infty)$.
The rest is straightforward.
\end{proof}

If we set $k=0$ in \eqref{relation 2} we obtain \eqref{relation 1}.

Let us denote the limit $q^\nu\rightarrow0$ symbolically by $\nu\rightarrow\infty$
(see Lemma \ref{rational}). From \eqref{relation 3} follows that
$$
\lim_{\nu\rightarrow\infty}J^{(k,\nu)}_{C,m}(q,z)=I_{C,m}(q,z|0,k).
$$
\subsection{Quasi-classical decomposition of fermionic sums}
In this subsection, we discuss the decomposition of fermionic sums with respect to
the dependence on some large parameters. We call it the quasi-classical decomposition.
We first explain the idea in simple examples, and then discuss more general cases.

\noindent{{\bf Example 1}.
Consider the simplest case $l=1,[r,s]=[0,k],\mu=0, m=1$. Then  one gets
$$I_{C,1}(q,z|0,k)=\frac{1+z+\cdots+z^k}{1-q}=\frac1{(1-q)(1-z)}-\frac{z^{k+1}}{(1-q)(1-z)}.$$
The result consists of two terms. It is a linear combination of $1$ and $z^k$
with rational function coefficients independent of $k$. This decomposition can be explained
by examining the large $k$ behaviour.
There is one particle in the interval $[0,k]$. The weight $p$ of the particle is restricted to
this interval, $0\leq p\leq k$.
If $k$ is large, the sum over $p$ for $0\leq p<\hskip-3pt<k$ and
that for $0\leq k-p<\hskip-3pt<k$ does not overlap.
Considering the sum over $p$ for $0\leq p$ and that for $p\leq k$ we obtain
$I_{C,1}(q,z|0,\infty)$ and $I_{C,1}(q,z|-\infty,k)$. In fact, the above decomposition is the same as
$$I_{C,1}(q,z|0,k)=I_{C,1}(q,z|0,\infty)+I_{C,1}(q,z|-\infty,k).$$
The sums in the right hand side contains more terms than the sum over  $0\leq p\leq k$.
However, since $\sum_{p\in\Z}z^p=0$ as rational function, we have the equality.

\noindent{{\bf Example 2}.
Let $l=1,[r,s]=[0,k],\mu=0$ as before but $m=2$. Since
$l=1$ the matrix $C$ is simply a scalar. We denote it by $c$. One can check
the following equality.
$$I_{C,2}(q,z|0,k)=I_{C,2}(q,z|0,\infty)+I_{C,1}(q,qz|0,\infty)I_{C,1}(q,z|-\infty,k)+I_{C,2}(q,z|-\infty,k).$$
Each of the terms in the right hand side has the distinction
that the dependence on $k$ enters only through $1,z^k,(z^2q^c)^k$, respectively.
These sums are over $p_1,p_2$ in the regions
\be
\{0\leq p_1\leq p_2\},\ \{0\leq p_1,0\leq k-p_2\},\ \{0\leq k-p_2\leq k-p_1\}.
\en
They are obtained by extending the following regions in the original sum:
\be
\{0\leq p_1\leq p_2<\hskip-3pt<k\},\ \{0\leq p_1<\hskip-3pt<k,0\leq k-p_2<\hskip-3pt<k\},\
\{0\leq k-p_2\leq k-p_1<\hskip-3pt<k\}.
\en

In general, we call such a decomposition quasi-classical.
We conjecture that the quasi-classical decompositions are exact for the fermionic sums.
In later sections, we prove the conjecture in some cases.

Here is the quasi-classical decomposition for the general case when the interval is finite.
\begin{conj}\label{I=IIconj}
We have
\begin{align}
I_{C,m}(q,z|r,s)=\left(z^mq^{W_{C,m}}\right)^r\sum_{0\le a\le m}\bigl(z^aq^{W_{C,a}}\bigr)^{s-r}
I_{C,m-a}(q,q^{Ca}z)I_{C,a}(q,z^{-1}q^{-Ca+{\rm diag}C}).\label{quasi classical}
\end{align}
\end{conj}

Let us explain how one can  obtain the right hand side.
Recall \eqref{quadratic} and \eqref{power z}.
Consider the fermionic sum on the interval $[r,s]$ where $r\rightarrow-\infty$ and $s\rightarrow\infty$.
 Let $0\leq a_i\leq m_i$ for $1\leq i\leq l$, and consider the sum over $p$ in the region
$r\leq p_{i,1}\leq\cdots\leq p_{i,m_i-a_i}<\hskip-3pt<s$ and
$r<\hskip-3pt<p_{i,m_i-a_i+1}\leq\cdots\leq p_{i,m_i}\leq s$.
We denote the vector of $p_{i,j}$ from the first region by ${\bf p}_{m-a}$
and from the second region by ${\bf p}_a$.
Two groups of variables are separated: if $p_{i,j}$ is in the first group and
$p_{i',j'}$ is in the second group, then we have
$$\min(p_{i,j},p_{i',j'})=p_{i,j}.$$
Therefore, we have
$$\overline Q_{C}({\bf p})=\overline Q_{C}({\bf p}_{m-a})+\overline Q_{C}({\bf p}_{a})
+\sum_{i'=1}^lC_{i,i'}a_{i'}\sum_{j=1}^{m_i-a_i}p_{i,j}.$$
Extending the regions for summation by removing the bounds of the form $\cdots<\hskip-3pt<s$ or
$r<\hskip-3pt<\cdots$, we obtain
\begin{align}
I_{C,m}(q,z|r,s)=\sum_{0\le a\le m}I_{C,m-a}(q,q^{Ca}z|r,\infty)I_{C,a}(q,z|-\infty,s).\label{two groups}
\end{align}
We rewrite this to \eqref{quasi classical} by using \eqref{J1} and \eqref{J2}.

In the special cases $r=s$ and $r=s+1$ Conjecture \ref{I=IIconj} reads as
follows.
\begin{align}
&\sum_{0\le a\le m}I_{C,m-a}(q,q^{Ca}z)I_{C,a}(q,z^{-1}q^{-Ca+{\rm diag}C})=\frac1{(q)_m},\\
&\sum_{0\le a\le m}\bigl(z^aq^{W_{C,a}}\bigr)^{-1}
I_{C,m-a}(q,q^{Ca}z)I_{C,a}(q,z^{-1}q^{-Ca+{\rm diag}C})=0.
\end{align}

The same quasi-classical decomposition procedure can be applied to $J^{(k,\nu)}_{C,m}(q,z)$,
where $k\rightarrow\infty$.
For this purpose we need fermionic sums on $\Z=(-\infty,\infty)$.
We define
\bea\label{pminfty}
X^{(k,\nu)}_{C,m}(q,z)=I_{C,\mu,m}(q,z|-\infty,\infty),
\ena
where $\mu$, which depends on $(k,\nu)$, is given by \eqref{corner}.

\begin{lem}
$X^{(k,\nu)}_{C,m}(q,z)$ is a rational function  in $q,z$.
It satisfies  the relations:
\be
X^{(k,\nu)}_{C,m}(q,z)=\left(z^mq^{W_{C,m}}\right)^kX^{(0,\nu)}_{C,m}(q,z)
\en
and
\bea
X^{(0,\nu)}_{C,m}(q,z)=X^{(0,\nu)}_{C,m}(q,z^{-1}q^{-\nu-Cm+{\rm diag}\,C}).\label{symmetry of X}
\ena
\end{lem}
\begin{proof}
The proof is straightforward. We only note that in order to prove
\eqref{symmetry of X} one needs to write the fermionic sum in
$\bf p$ variables and  change the summation variable $p_{i,j}$ to
$-p_{i,m_i+1-j}$ for all $i,j$.
\end{proof}

The following proposition is an analogue of Lemma \ref{J=II}.
\begin{prop}
\begin{align}
X^{(0,\nu)}_{C,m}(q,z)&=\sum_{0\le a\le m}\left(z^aq^{W_{C,a}+Ca\cdot(m-a)}\right)^{-1}
I_{C,a}(q,z^{-1}q^{-Cm+{\rm diag}C})I_{C,m-a}(q,q^\nu z),\label{X1}\\
&=\sum_{0\le a\le m}z^aq^{W_{C.a}+\nu\cdot a}I_{C,m-a}(q,z^{-1}q^{-Cm+{\rm diag}C})I_{C,a}(q,q^\nu z).\label{X2}
\end{align}
\end{prop}
\begin{proof}
We use two cuttings of the infinite interval of weights $(-\infty,\infty)$
into semi-infinite intervals:
$(-\infty,\infty)=(-\infty,-1]\sqcup[0,\infty)$ or
$(-\infty,\infty)=(-\infty,0]\sqcup[1,\infty)$.
This leads to \eqref{X1} and \eqref{X2}.
\end{proof}

Applying the procedure of the quasi-classical decomposition, we obtain
\begin{conj}
\be
J^{(k,\nu)}_{C,m}(q,z)
=\sum_{0\le a\le m}I_{C,m-a}(q,q^{Ca}z)X^{(k,\nu)}_{C,a}(q,z).
\en
\end{conj}
In the right hand side the summation variable $a_i$ where $a=(a_1,\ldots,a_l)$,
represents the number of color $i$ particles whose weights are "close to $k$".
The weights of the remaining particles are "small" compared to $k$.
We conjecture that this is exact for finite $k$.
In particular, setting $k=0$ and $k=-1$ and using \eqref{relation 1}, we obtain
\bea
I_{C,m}(q,q^\nu z)
&=&\sum_{0\le a\le m}I_{C,m-a}(q,q^{Ca}z)X^{(0,\nu)}_{C,a}(q,z),\label{recursion J=JX}\\
q^{\nu\cdot m}I_{C,m}(q,q^\nu z)
&=&\sum_{0\le a\le m}I_{C,m-a}(q,q^{Ca}z)\left(z^aq^{W_{C,a}}\right)^{-1}
X^{(0,\nu)}_{C,a}(q,z).\label{recursion J=JX k=-1}
\ena
In Section \ref{sec:4} we prove these equalities in the case where $C$
is a simply-laced Cartan matrix. We also give a generalization of these equalities
in the case where $C$ is the symmetrization of a non simply-laced Cartan matrix.

Finally, we give the quasi-classical decomposition for the fermionic sum on $(-\infty,\infty)$
with two corners at $0$ and $k$, with angle $\nu_1$ and $\nu_2$, respectively.

We denote this quantity by $X^{(0,\nu_1;k,\nu_2)}_{C,m}(q,z)$. We conjecture that
$$
X^{(0,\nu_1;k,\nu_2)}_{C,m}(q,z)=
\sum_{0\le a\le m}\left(z^aq^{W_{C,a}}\right)^kX^{(0,\nu_1)}_{C,m-a}(q,q^{Ca}z)X^{(0,\nu_2)}_{C,a}(q,z).
$$
Restricting to $k=0$ we have
\bea\label{X=XX}
X^{(0,\nu_1+\nu_2)}_{C,m}(q,z)=
\sum_{0\le a\le m}X^{(0,\nu_1)}_{C,m-a}(q,q^{Ca}z)X^{(0,\nu_2)}_{C,a}(q,z)
\ena

\subsection{The case of $\frak{sl}_2$}
In this subsection, we restrict to the $\frak{sl}_2$ case, i.e., $l=1$ and $C=2$,
and write some of the fermionic sums and their relations explicitly.
Moreover, we discuss vanishing theorems which arise in connection with integrality
of angle variables. In the following we drop $C$ in the notation because it is fixed to $C=2$.
We also use the $q$ binomial coefficients defined by
\be
\left[\gamma\atop n\right]=\frac{(q^{\gamma-n+1})_n}{(q)_n}.
\en
Here $n$ is a non-negative integer, but $\gamma$ is arbitrary, possibly a
formal variable.

First, we recall a known result and its proof (see \cite{FFJMM}).
\begin{prop}\label{prop 1}
We have
\begin{align}
I_m(q,z)=\frac1{(q)_m(z)_m}.\label{sl2 J}
\end{align}
\end{prop}
\begin{proof}
The recursion \eqref{fermionic recursion} in this case reads as
$$I_m(q,z)=\sum_{0\le a\le m}\frac{z^aq^{a(a-1)}}{(q)_{m-a}}I_a(q,z).$$
The fermionic sum is uniquely determined by this recursion with the initial condition $I_0(q,z)=1$.
Therefore, it is enough to prove this recursion for \eqref{sl2 J}.
After substitution, we want to prove
$$\frac1{(z)_m}=\sum_{0\le a\le m}\left[m\atop a\right]\frac{z^aq^{a(a-1)}}{(z)_a}.$$
Using
$\left[m\atop a\right]=q^a\left[m-1\atop a\right]+\left[m-1\atop a-1\right]$, we obtain
\begin{align*}
(RHS)&=\sum_{a=0}^{m-1}{\textstyle\left[m-1\atop a\right]}
\left(\frac{z^aq^{a^2}}{(z)_a}+\frac{z^{a+1}q^{a(a+1)}}{(z)_{a+1}}\right)\\
&=\frac1{1-z}\sum_{a=0}^{m-1}
\frac{\left[m-1\atop a\right](qz)^aq^{a(a-1)}}{(qz)_a}=(LHS).
\end{align*}
\end{proof}
The following proposition holds for an arbitrary value of $\nu$.
\begin{prop}\label{prop 2}
$$
X^{(0,\nu)}_m(q,z)=\frac{\left[\nu\atop m\right]}{(z^{-1}q^{2(1-m)})_m(q^\nu z)_m}.
$$
\end{prop}
\begin{proof}
We use the decomposition \eqref{X2}. Substituting the expressions for
$I_{C,a}(q,z|-\infty,k)$ and $I_{C,m-a}(q,q^\nu z)$ given by \eqref{J2} and \eqref{sl2 J}, we obtain
$$X^{(0,\nu)}_m(q,z)=\frac1{(q)_m}\sum_{0\le a\le m}\frac{\left[m\atop a\right](q^\nu z)^aq^{a(a-1)}}
{(z^{-1}q^{2-2m})_{m-a}(q^\nu z)_a}.$$
The rest of the proof goes similarly as in Proposition \ref{prop 1}.
\end{proof}
If $n$ is a non-negative integer, the range for $m$ where $X^{(0,n)}_m(q,z)$ is non-zero is restricted.
\begin{cor}
If $n\in\Z_{\geq0}$ and $m>n$, then we have $X^{(0,n)}_m(q,z)=0$.
\end{cor}
We identify $n$ with the highest weight of the irreducible representation $V_n$ of $\frak{sl}_2$.
The above statement says the fermionic sum
is non-zero only if $n-2m$ is a weight of $V_n$. In Proposition \ref{lem:vanish}
we establish
vanishing theorems of the form $X^{(0,\nu)}_{C,m}=0$
for the case where $C$ is a simply-laced Cartan matrix, the parameter $\nu$ corresponds
to a dominant integral weight, and $\nu-\sum_im_i\alpha_i$ is not a weight of $V_\nu$. Our tool is
the representation theory of $U_v(\g)$ with $v^2=q$, where $\g$ is a simple Lie
algebra associated with $C$. We expect vanishing theorems of this kind
are valid in a much wider class of $C$, though it is beyond the scope of this paper.

Two recursions \eqref{recursion J=JX} and \eqref{recursion J=JX k=-1} reads as follows.
\begin{prop} \label{prop 3}
For an arbitrary value of $\nu$, we have
\begin{align}
I_m(q,q^\nu z)=\sum_{0\le a\le m}\frac{\left[\nu\atop a\right]}{(z^{-1}q^{2(1-a)})_a(q^\nu z)_a}I_{m-a}(q,q^{2a}z),\\
q^{\nu m}I_m(q,q^\nu z)=\sum_{0\le a\le m}\frac{z^{-a}q^{-a(a-1)}\left[\nu\atop a\right]}
{(z^{-1}q^{2(1-a)})_a(q^\nu z)_a}I_{m-a}(q,q^{2a}z).
\end{align}
\end{prop}

The quasi-classical decomposition generates many more identities than we discussed above.
Here we give an example. Let $0\leq k_1\leq\cdots\leq k_n$ be non-negative integers.
Set
$$\mu_t=\sum_{i=1}^n(t-k_i)_+\quad
\hbox{where }(t)_+=\begin{cases}t&\hbox{ if }t\geq0;\\0&\hbox{ otherwise.}\end{cases}$$
Namely, the linear coefficient $\mu$ has corners at $k_i$ with angle $\sum_{r=1}^n\delta_{k_i,k_r}$.
Suppose that the quasi-classical decomposition is exact. Then for $n\ge m$ we have
\begin{align}
I_{\mu,m}(q,z)&=\sum_{0\le a\le m}I_a(q,q^{2(m-a)}z)X_{\mu,m-a}(q,z),\\
X_{\mu,m}(q,z)&=\sum_{\varepsilon_1,\ldots,\varepsilon_n=0,1\atop\varepsilon_1+\cdots+\varepsilon_n=m}
X^{k_1,\ldots,k_n}_{\varepsilon_1,\ldots,\varepsilon_n},\\
X^{k_1,\ldots,k_n}_{\varepsilon_1,\ldots,\varepsilon_n}&=\prod_{i:\ \varepsilon_i=1}
\frac{q^{-\sum_{r=1}^{i-1}k_r}(q^{\varepsilon(i)}z)^{k_i}}{g(q^{\varepsilon(i)}z)}.
\end{align}
Here $g(z)=(1-z^{-1})(1-qz)$, $\varepsilon(i)=2\sum_{r=i+1}^n\varepsilon_r+i-1$.
Note that $X^{(0,1)}_{1,0}=1/g(z)$.

For example, setting $n=3,k_1=k_2=k_3=0,m=2$, we obtain
$$\frac{\left[3\atop 2\right]}{(q^3z)_2(q^{-2}z^{-1})_2}=
\frac1{g(zq^2)}\left(\frac1{g(zq)}+\frac1{g(zq^2)}+\frac1{g(zq^3)}\right).$$

%%%%%%%%%%%%%%%%%%%%%%%%%%%%%%%%%%%%%%

\section{Fermionic and Toda recursions}\label{sec:3}

In this section, we develop the representation theoretic
approach to fermionic formulas.
In certain cases we show that they coincide with scalar products of
Whittaker vectors, which are
eigenfunctions of the quantum Toda Hamiltonian.

Quantum deformation of Whittaker vectors has been
introduced and studied by Sevostyanov \cite{Sev}. An independent
construction was given by Etingof \cite{Et}
from a slightly different point of view.
A geometric interpretation of the
eigenvectors of the quantum Toda Hamiltonian
as Shapovalov scalar product of Whittaker vectors
has been given by Braverman and Finkelberg
\cite{BrFi}, reproducing the
main results of Givental and Lee \cite{GiL}.
In this section we give a review of this subject,
following closely the algebraic framework of \cite{Sev}
with minor modifications.
We shall show that fermionic formulas naturally
arise from Drinfeld's Casimir element.

\subsection{Quantum groups}

We fix the notation as follows.
Let $\g$ be a complex simple Lie algebra,
$\h$ the Cartan subalgebra, $\al_1,\cdots,\al_l$
the simple roots and $\omega_1,\cdots,\omega_l$
the fundamental weights. Set
$P=\oplus_{i=1}^l\Z\omega_i$,
$Q=\oplus_{i=1}^l\Z\al_i$,
$P_+=\oplus_{i=1}^l\Z_{\ge0}\omega_i$,
$Q_+=\oplus_{i=1}^l\Z_{\ge0}\al_i$.
Let further
$\Delta_+$ denote the set of positive roots.
We fix
a non-degenerate invariant bilinear form
$(~,~):\h\times\h\to\C$
such that $(P,Q)\subset\Z$,
and identify $\h^*$ with $\h$ via $(~,~)$.
We set
\be
d_i=\frac{1}{2}(\al_i,\al_i), \quad
\al_i^{\vee}=d_i^{-1}\al_i,
\quad
\rho=\sum_{r=1}^l\omega_r.
\en
We choose $(~,~)$ so that $d_1,\cdots,d_l$ are
relatively prime positive integers.

Let $\mathcal{N}$
be a positive integer satisfying
$(P,P)\subset (1/\mathcal{N})\Z$.
The quantum group $\Uv$ is a unital associative
algebra over $\K=\C(v^{1/\mathcal{N}})$, with generators
$E_i,F_i$ ($1\le i\le l$),
$K_\mu$ ($\mu\in P$) and the standard defining relations
\be
&&K_\mu K_{\mu'}=K_{\mu+\mu'},\quad K_0=1, \\
&&K_\mu E_i K_\mu^{-1}=v^{(\mu,\al_i)}E_i,
\quad K_\mu F_i K_\mu^{-1}=v^{-(\mu,\al_i)}F_i,
\\
&&[E_i,F_j]=\delta_{ij}\frac{K_i-K_i^{-1}}{v_i-v_i^{-1}},
\\
&&\sum_{s=0}^r(-1)^sE_i^{(r-s)}E_jE_i^{(s)}=0
\qquad (r=1-(\al^{\vee}_i,\al_j),\ i\neq j),
\\
&&\sum_{s=0}^r(-1)^s F_i^{(r-s)}F_jF_i^{(s)}=0
\qquad (r=1-(\al^{\vee}_i,\al_j),\ i\neq j).
\en
Here $K_i=K_{\al_i}$, $v_i=v^{d_i}$, $X_i^{(s)}=X_i^s/[s]_{v_i}!$
($X=E,F$) and  $[s]_v!=\prod_{p=1}^s(v^p-v^{-p})/(v-v^{-1})$.
We choose the coproduct
\be
&&\Delta E_i=E_i\otimes 1+K_i\otimes E_i,
\\
&&\Delta F_i=F_i\otimes K_i^{-1}+1\otimes F_i,
\\
&&\Delta K_\mu=K_\mu\otimes K_\mu,
\en
antipode
\be
S(E_i)=-K_i^{-1}E_i,\ S(F_i)=-F_i K_i,\ S(K_\mu)=K_\mu^{-1}
\en
and counit
\be
\varepsilon(E_i)=\varepsilon(F_i)=0,\
\varepsilon(K_\mu)=1.
\en

We shall also consider the quantum group $\Uvv$ with
parameter $v^{-1}$. Denote
the generators by $\bar E_i$, $\bar F_i$, $\bar K_\mu$.
We choose the opposite coproduct,
\be
&&\Delta \bar E_i=\bar E_i\otimes \bar K_i+1 \otimes \bar E_i,
\\
&&\Delta \bar F_i=\bar F_i\otimes 1+
\bar K_i^{-1}\otimes \bar F_i,
\\
&&\Delta \bar K_\mu=\bar K_\mu\otimes \bar K_\mu,
\en
antipode
\be
S(\bar E_i)=-\bar E_i\bar K_i^{-1},
\ S(\bar F_i)=-\bar K_i\bar F_i ,\
S(\bar K_\mu)=\bar K_\mu^{-1},
\en
and counit
\be
\varepsilon(\bar E_i)=\varepsilon(\bar F_i)=0,\
\varepsilon(\bar K_\mu)=1.
\en

There is a $\K$-linear anti-isomorphism of algebras
given by
\bea
\sigma:\Uv\to \Uvv,
\quad E_i\mapsto \bar F_i,\
F_i\mapsto \bar E_i,
\ K_\mu \mapsto \bar K_\mu^{-1}.
\label{sigma}
\ena
We have
\bea
\Delta\circ\sigma=\sigma\otimes\sigma\circ\Delta\,.
\ena

\subsection{Verma modules}

For $\la\in P$, let $\V^\la$
be the Verma module over $\Uv$ generated by the
highest weight vector $\one^\la$ with defining relations
\be
E_i\one^\la=0\quad (1\le i\le l),
\qquad K_\mu\one^\la=v^{(\mu,\la)}\one^\la\quad (\mu\in P).
\en
Similarly let $\bV^\la$
be the Verma module over $\Uvv$ generated by the
highest weight vector $\bar\one^\la$ with defining relations
\be
\bar E_i\bar \one^\la=0\quad (1\le i\le l),
\qquad \bar K_\mu\bar \one^\la=v^{-(\mu,\la)}\bar
\one^\la\quad (\mu\in P).
\en
We have obvious gradings
$\V^\la=\oplus_{\beta\in Q^+}(\V^\la)_\beta$,
$\bV^\la=\oplus_{\beta\in Q^+}
(\bV^\la)_\beta$,
where %$Q_+=\sum_{r=1}^l\Z_{\ge 0}\al_r$ and
\bea
&&(\V^\la)_{\beta}=\{w\in \V^\la\mid
K_\mu w=v^{(\mu,\la-\beta)}w\quad(\mu\in P)
\},
\label{wt1}\\
&&(\bV^\la)_{\beta}=
\{w\in \bV^\la\mid
\bar K_\mu w=v^{-(\mu,\la-\beta)}w\quad(\mu\in P)\}.
\label{wt2}
\ena
There exists a unique non-degenerate $\K$-bilinear pairing
$(~,~):\V^\la\times\bV^\la\to\K$, such that
$(\one^\la,\bar\one^\la)=1$ and
\bea
(x w, w')=(w,\sigma(x)w')
\label{Shapo}
\ena
for all $x\in\Uv$ and $w\in\V^\la$,
$w'\in \bV^\la$.
The weight components \eqref{wt1}, \eqref{wt2}
are mutually orthogonal with respect to $(~,~)$.
We extend the scalar product on tensor products of Verma
modules as
\be
(u_1\otimes u_2,v_1\otimes v_2)=(u_1,v_1)(u_2,v_2)
\quad (u_i\in \V^{\lambda_i},\ v_i\in \bar\V^{\lambda_i})\,.
\en

\subsection{Whittaker vectors}

Whittaker vectors are defined by giving the following data:
an orientation of the Dynkin graph,
a set of elements $\nu_i\in P$,
and non-zero scalars $c_i\in\K$.
An orientation is represented by
a skew-symmetric matrix $\e=(\e_{i,j})$, where
$\e_{i,j}=1$ if there is an arrow pointing
from node $i$ to node $j$,
and $\e_{i,j}=0$ if $i$ and $j$ are disconnected.
Given $\e$, choose $\nu_i\in P$ such that
\bea
(\nu_i,\al_j)-(\nu_j,\al_i)=\e_{i,j}(\al_i,\al_j).
\label{nij}
\ena
For instance one can take
$\nu_i=\sum_{k=1}^{i-1}\omega_k
\e_{i,k}(\al_i,\al^{\vee}_k)$.

A Whittaker vector associated with the data
$\e$, $\nu=(\nu_i)$ and $c=(c_i)$ is an element
\be
\theta^\la=\theta^\la(\e,\nu,c)=\sum_{\beta\in Q_+}
\theta^\la_\beta,
\qquad \theta^\la_\beta\in (\V^\la)_\beta,
\en
which belongs to a completion
$\widehat{\V}^\la=\prod_{\beta\in Q_+}(\V^\la)_\beta$
of the Verma module.
It is defined by the conditions
$\theta^\la_0=\one^\la$ and
\bea
E_i K_{\nu_i}\ \theta^\la
=\frac{c_i}{1-v_i^2}\ \theta^\la.
\label{Whit}
\ena
It is known that the Whittaker vector exists and is unique
\cite{Sev}.
Fixing $\e$ and changing $\nu,c$ results in
a simple scalar multiple of the weight components of
$\theta^\la(\e,\nu,c)$.
Explicitly we have the transformation law
\bea
&&\theta^\la_\beta(\e,\nu+\ga,c')=
v^{(1/2)\sum_{i=1}^l(\beta,\omega^{\vee}_i)(\beta-2\lambda,\ga_i)}
\theta^{\la}_\beta(\e,\nu,c),
\label{trlaw1}\\
&&\theta^\la_\beta(\e,\nu,c'')
=v^{-(\kappa,\beta)}\theta^\la_\beta(\e,\nu,c),
\label{trlaw2}
\ena
where
$\ga_i\in P$ is such that
$(\ga_i,\al_j)=(\ga_j,\al_i)$ holds, and
$c'_i=v^{-(\ga_i,\al_i)/2}c_i$, $c''_i=v^{-(\kappa,\al_i)}c_i$.

Similarly one defines the dual Whittaker vector
\be
&&\bar\theta^\la=\bar\theta^\la(\e,\nu,c)
%\in \prod_\beta(\bV^\la)_\beta
\in \widehat{\bar \V}^\la,
\qquad
\widehat{\bar \V}^\la=\prod_{\beta\in Q_+}(\bV^\la)_\beta,
\en
imposing $\bar\theta^\la_0=\bar\one^\la$ and
\bea
\bar E_i \bar K_{\nu_i}\ \bar \theta^\la
=\frac{c_i^{-1}}{1-v_i^{-2}}\ \bar \theta^\la
\label{DWhit}
\ena
in place of \eqref{Whit}.

The main object of our interest is the scalar product
\bea
J^\la_\beta=v^{-(\beta,\beta)/2+(\la,\beta)}\
(\theta^\la_\beta(\e,\nu,c),\bar\theta^\la_\beta(\e,\nu,c)).
\label{Jla}
\ena
We set $J^\la_\beta=0$ unless $\beta\in Q_+$.

From \eqref{trlaw1} and \eqref{trlaw2}, we see that
\eqref{Jla} is independent of the choice of $\nu, c$.
As it will turn out, it is actually independent of
the orientation $\e$ as well (see Theorem \ref{Frec}
below). Anticipating this fact, we suppress the dependence
on $\e,\nu,c$ from the notation.

\subsection{Fermionic recursion}

Now we state the main result of this section.
Set $q=v^2$, $q_i=q^{d_i}$, and
\be
&&(q)_\beta=\prod_{i=1}^l(q_i;q_i)_{m_i}\quad
\mbox{ for $\beta=\sum_{i=1}^lm_i\al_i$ }.
\en

The following is a counterpart of the
fermionic recursion given in Proposition \ref{sec2:frec}.
\begin{thm}\label{Frec}
The quantities $J^\la_\beta$ are
uniquely characterized by $J^\la_0=1$
and the recursion relation
\bea
J^\la_\beta=\sum_{\ga}
\frac{1}{(q)_{\beta-\ga}}
q^{(\ga,\ga)/2-(\la+\rho,\ga)}J^{\la}_\ga.
\label{eq:frec}
\ena
In particular, $J^\la_\beta$ is independent of the choice of
orientation which enters the definition.
\end{thm}
The $J^\la_\beta$ are determined as
rational functions in the variables
$q$ and $z_i=q^{-(\la,\al_i)}$.

In accordance with the previous section,
let us introduce the following sum
for a (possibly infinite) interval $[r,s]$.
\bea
J^\la_\beta[r,s]=
\sum_{\sum_{t=r}^s\ga^{(t)}=\beta}
\frac{q^{B(\ga)}}{\prod_{t=r}^s(q)_{\ga^{(t)}}},
\label{eq:frec2}
\ena
where
\be
&&
B(\ga)=
\frac{1}{2}\sum_{r\le t,t'\le s}\min(t,t')
(\ga^{(t)},\ga^{(t')})
-(\la+\rho,\sum_{t=r}^s t \ga^{(t)}).
\en
We have
\bea
&&J^\la_\beta[0,0]=\frac{1}{(q)_\beta},
\label{Jrs1}\\
&&J^\la_\beta[r+1,s+1]=q^{(\beta,\beta)/2-(\la+\rho,\beta)}J^\la_\beta[r,s],
\label{Jrs2}\\
&&J^\la_\beta[r,s]=\sum_{\al+\ga=\beta}J^{\la-\ga}_\al[r,u]J^\la_\ga[u+1,s]
\qquad (r\le u<s).
\label{Jrs3}
\ena

\begin{thm}\label{Frec2}
In the notation above, we have
\be
J^\la_\beta=J^\la_\beta[0,\infty).
\en
\end{thm}
\begin{proof}
This is a restatement of Theorem \ref{Frec}.
\end{proof}

Suppose $C$ is a simply laced Cartan matrix,
and let $\beta=\sum_{i=1}^lm_i\al_i$, $z_i=q^{-(\la,\al_i)}$.
Since $(Cm,m)=(\beta,\beta)$ and
$({\rm diag} C,m)=(2\rho,\beta)$, we have
\be
J^\la_\beta[0,\infty)=I_{C,m}(q,z)
\en
where the right hand side is defined in \eqref{ICm}.
We shall discuss the interpretation of
$J^\la_\beta[r,s]$ for finite interval $[r,s]$
in the next section (see Theorem \ref{thm:J[0,k]}).
When $C$ is non-simply laced, \eqref{eq:frec2} gives
a generalization of the fermionic sum considered
in the previous section
due to the denominator $(q_i;q_i)_n$.

\begin{cor}\label{propJ}
$(\hbox{\rm i})$
The rational function $J^\la_\beta$
is regular at $z_1=\cdots=z_l=0$ and
\bea
{J^\la_\beta}\bigl|_{z_1=\cdots=z_l=0}=
\prod_{r=1}^l
\frac{1}{(q_i;q_i)_{m_i}}.
\label{Jz=0}
\ena

$(\hbox{\rm ii})$
We have the symmetry property
\bea
J^\la_\beta\Bigl|_{q\to q^{-1}}=
q^{(\beta,\beta)/2-(\la,\beta)}J^\la_\beta.
\label{Jqinv}
\ena

$(\hbox{\rm iii})$
The set $\{J^\la_{\beta}\}_{\beta\in Q_+}$
is linearly independent over $\C(q)$.
\end{cor}
\begin{proof}
Assertion (i) is a direct consequence of \eqref{eq:frec2} since $J^\la_\beta=J^\la_\beta[0,\infty)$.
In the definition of $J^\la_\beta$,
$\theta^\la$ and $\bar\theta^\la$ enter in a symmetric way.
Therefore assertion (ii) follows from the definition \eqref{Jla}.
(Note that the change of variable $q\to q^{-1}$ in the left hand side of
(ii) implies $z_i\to z_i^{-1}$).

To see (iii), it suffices to show that
$J^\la_{\beta}\Bigl|_{q\to q^{-1}}$ constitute a
linearly independent set.
Property (ii) implies that
each of them has distinct leading power
$z_1^{m_1}\cdots z_l^{m_l}$ in $z_1,\cdots,z_l$.
Hence the independence is evident.
\end{proof}

\subsection{Derivation of the fermionic recursion}

In this subsection we give a derivation of Theorem
\ref{Frec}.

First we recall the Cartan-Weyl basis and the
product formula for the universal $R$ matrix due to
Khoroshkin and Tolstoy \cite{KT}.
By definition, a total order $<$ on $\Delta_+$ is said to be
normal if
$\al,\beta,\al+\beta\in \Delta_+$ and $\al<\beta$
imply $\al<\al+\beta<\beta$.
Normal orders are in one-to-one correspondence with
reduced decompositions of the longest element of the Weyl group.
Moreover an arbitrary total order on the set of simple roots
can be extended to a normal order on $\Delta_+$
\cite{Tol,Zh}.
To a normal order $<$, one associates root vectors
\bea
E^{<}_\beta,\ F^{<}_\beta \qquad (\beta\in\Delta_+)
\label{root-vec}
\ena
by induction on $h(\beta)$, where
$h(\sum_{i=1}^ln_i\al_i)=\sum_{i=1}^l n_i$.
When $h(\beta)=1$,
define $E^{<}_{\al_i}=E_i$,  $F^{<}_{\al_i}=F_i$.
Let $\ga$ be an element with
$h(\ga)=n$, and suppose that
\eqref{root-vec} are already defined for $h(\beta)<n$.
Choose a decomposition $\ga=\al+\beta$ in such a way that
there are no other roots $\al',\beta'\in\Delta_+$ satisfying
$\ga=\al'+\beta'$, $\al\le \al'<\beta'\le \beta$. Then define
\bea
&&E^{<}_\ga=E^{<}_\al E^{<}_\beta
-v^{(\al,\beta)}E^{<}_\beta E^{<}_\al,
\label{root-vec2}
\\
&&F^{<}_\ga=
c_\ga(F^{<}_\beta F^{<}_\al
-v^{-(\al,\beta)}F^{<}_\al F^{<}_\beta).
\label{root-vec3}
\ena
The scalar $c_\ga\in \K$ can be so chosen that
$[E^{<}_\ga,F^{<}_\ga]= (K_\ga-K_\ga^{-1})/(v_\ga-v_\ga^{-1})$.
Here and after we set $v_\ga=v^{(\ga,\ga)/2}$.

The product formula is written in terms of the $q$-exponential
function
\be
&&\exp_q(x)=\sum_{n=0}^\infty \frac{((1-q)x)^n}{(q;q)_n},
\en
which satisfies $\exp_q(x)\exp_{q^{-1}}(-x)=1$.
\begin{prop}\cite{KT}\label{R}
Fix a normal order $<$, and consider the element
\bea
&&\Theta=
\overset{\curvearrowleft}{\prod_{\beta\in\Delta_+}
}\exp_{v_\beta^{2}}
\left(-(v_\beta-v_\beta^{-1})
F^{<}_\beta\otimes E^{<}_\beta\right),
\label{Rprod}
\ena
where the product is so ordered that $\beta$ appears to
the right of $\beta'$ if $\beta<\beta'$.
Then $\Theta$
does not depend on the choice of the normal order $<$.
\end{prop}
The universal $R$ matrix is given by
$\R=\Theta_{21}v^{-T}$, where
$\Theta_{21}=\tau(\Theta)$, $\tau(a\otimes b)=b\otimes a$,
and $T\in \h\otimes\h$ stands for the canonical element.

The following construction is well known.
\begin{prop}\label{DCas}\cite{Dri}
Set $u=m(S\otimes \id)\Theta$, where $m(a\otimes b)=ab$
stands for the multiplication.
Then $u$ is a well defined operator on $\V^\la$.
It acts on each $(\V^\la)_\beta$ as a scalar,
\be
u\bigl|_{(\V^\la)_\beta}=v^{-(\beta,\beta)+2(\la+\rho,\beta)}
\times\id_{(\V^\la)_\beta}.
\en
\end{prop}
The formal element
$v^{\sum_{i=1}^l\al_i^{\vee}\omega_i+2\rho}u$ is
sometimes referred to as the Drinfeld (quantum) Casimir element.
It acts on $\V^\la$ as a scalar: $v^{\sum_{i=1}^l\al_i^{\vee}\omega_i+2\rho}u\bigl|_{\V^\la}
=v^{(\la+2\rho,\la)}\times\id_{\V^\la}$.

A normal order $<$ is said to be compatible with an orientation
$\e$ of the Dynkin graph
if $\e_{j,i}=1$ implies $\al_i<\al_j$.
\begin{prop}\label{Simple}\cite{Sev}
Let \eqref{root-vec} be the root vectors with respect to
a normal order $<$, and let $\theta^\la(\e,\nu,c)$ be
the Whittaker vector.
If the order $<$ is compatible with the orientation $\e$, then
we have
\bea
E^{<}_\beta \ \theta^\la(\e,\nu,c)=0
\qquad \mbox{for all non-simple roots $\beta\in\Delta_+$}.
\label{root-vanish}
\ena
\end{prop}
\begin{proof}

Since the root vectors are defined
by multiple $v$-commutators \eqref{root-vec2},
\eqref{root-vec3}, it is sufficient to show
\eqref{root-vanish} when $\ga=\al_i+\al_j$. This can be
verified by a direct calculation using \eqref{nij} and
\eqref{Whit}.
\end{proof}

\noindent{\it Proof of Theorem \ref{Frec}.}\quad
By Proposition \ref{DCas}, we have
\be
(u\ \theta^\la_\beta,\bar \theta^\la_\beta)
=v^{-(\beta,\beta)+2(\la+\rho,\beta)}
(\theta^\la_\beta,\bar \theta^\la_\beta).
\en
Suppose that $\ga_1<\cdots<\ga_l$ are
the simple roots appearing in the chosen normal order.
Expanding the formula \eqref{Rprod}, we obtain
\be
&&u\ \theta^\la_\beta=
\sum_{n_1,\cdots,n_l}
\prod_{i=1}^l
\left(
\frac{v_i^{-n_i}(1-v^2_i)^{2n_i}}{(v_i^2;v_i^2)_{n_i}}\right)
\\
&&
\times
S(F_{\ga_1})^{n_1}\cdots S(F_{\ga_l})^{n_l}
E_{\ga_l}^{n_l}\cdots E_{\ga_1}^{n_1}\ \theta^\la_\beta.
\en
In view of Proposition \ref{Simple},
we have retained only those terms consisting of
simple roots.
Setting $\ga=\sum_{i=1}^ln_i\ga_i$ and using \eqref{Shapo},
we obtain
\begin{align*}
&(u\ \theta^\la_\beta,\bar \theta^\la_\beta)
=
\sum_{n_1,\cdots,n_l}
\prod_{i=1}^l
\left(\frac{(-v_i)^{-n_i}(1-v^2_i)^{2n_i}}{(v_i^2;v_i^2)_{n_i}}
\right)
v^{(\ga,\ga)/2+(\ga,\la+\rho-\beta)}
\times
(E_{\ga_l}^{n_l}\cdots E_{\ga_1}^{n_1}\ \theta^\la_\beta,
\bar E_{\ga_l}^{n_l}\cdots
\bar E_{\ga_1}^{n_1}\bar\theta^\la_\beta).
\end{align*}

%
%
%\be
%(u\ \theta^\la_\beta,\bar \theta^\la_\beta)
%=
%\sum_{n_1,\cdots,n_l}
%\prod_{i=1}^l
%\left(\frac{(-v_i)^{-n_i}(1-v^2_i)^{2n_i}}{(v_i^2;v_i^2)_{n_i}}
%\right)
%v^{(\ga,\ga)/2+(\ga,\la-\beta-\rho)}
%\times
%(E_{\ga_l}^{n_l}\cdots E_{\ga_1}^{n_1}\ \theta^\la_\beta,
%\bar E_{\ga_l}^{n_l}\cdots
%\bar E_{\ga_1}^{n_1}\bar\theta^\la_\beta).
%\en
Now apply the relations following from
\eqref{Whit}, \eqref{DWhit},
\be
v^{(\nu_i,\la-\beta)}E_i\theta^\la_\beta=
\frac{c_i}{1-v_i^2}\ \theta^\la_{\beta-\al_i},
\quad
v^{-(\nu_i,\la-\beta)}\bar E_i\bar \theta^\la_\beta=
\frac{c_i^{-1}}{1-v_i^{-2}}\ \bar \theta^\la_{\beta-\al_i}.
\en
Since the generators $E_{\ga_i}$, $\bar E_{\ga_i}$ are
arranged in the same order,  the powers of $v$ cancel out.
Substituting the definition \eqref{Jla} of $J^\la_\beta$
and simplifying the formula, we obtain the desired result.
\qed

\subsection{Toda recursion}
In \cite{Sev,Et}, 
quantum difference Toda Hamiltonians have been derived from 
Whittaker functions on quantum groups. 
For example, when $\g=\mathfrak{sl}_{l+1}$, the simplest 
Hamiltonian is given by (see \cite{Et}, (6.5))
\begin{align*}
&H_{Toda}=\sum_{i=1}^{l+1}D_{i}-(q-q^{-1})^2 \sum_{i=1}^l y_i^{-1}y_{i+1}D_{i},
%H_{Toda}=\sum_{i=0}^lD_{i+1}D_i^{-1}(1-y_i)\prod_{j=i+1}^l\left(q^{-1}z_j\right),
%\quad
%\varepsilon=\sum_{i=0}^l\prod_{j=i+1}^l\left(q^{-1}z_j\right),
%\label{TodaHamiltonian}
\end{align*}
where $D_i$ denotes the shift operator
$(D_i f)(y_1,\ldots,y_{l+1})=f(y_1,\ldots,qy_i,\ldots,y_{l+1})$.
%with $D_0=D_{l+1}=1,\,y_0=y_{l+1}=0$.
In this subsection we show that the generating function
\begin{align*}
%F(q,z,y)=\sum_{\beta}J^\la_\beta\prod_{i=1}^ly_i^{(\beta,\omega_i)/d_i}
F(q,z,y)=\sum_{\beta}\bigl(q(q-q^{-1})^2\bigr)^{|\beta|}
\prod_{i=1}^{l+1}y_i^{(\lambda-\beta,\epsilon_i)-i}\, J^\la_\beta 
\end{align*}
is an eigenfunction of $H_{Toda}$
%the quantum difference Toda Hamiltonian:
%$H_{Toda}F=\varepsilon F$ 
with eigenvalue $\varepsilon=\sum_{i=1}^{l+1}q^{(\lambda,\epsilon_i)-i}$.
Here $\alpha_i=\epsilon_i-\epsilon_{i+1}$ and $|\beta|=\sum_{i=1}^ln_i$ if 
$\beta=\sum_{i=1}^ln_i\alpha_i$. 
More specifically we show this in the form of 
an equivalent set of recurrence relations for $J^\la_\beta$;
see formulas below Proposition \ref{TodaRec}.

The argument to derive these relations is similar to the one for Theorem \ref{Frec}.
In place of the Drinfeld Casimir element, we use
other central elements including the quadratic Casimir element.

%
%As shown in \cite{Sev,Et},
%the generating function
%\be
%F(q,z,y)=\sum_{\beta}J^\la_\beta\prod_{i=1}^l
%y_i^{(\beta,\omega_i)/d_i}
%\en
%is an eigenfunction of the
%quantum difference Toda Hamiltonian:
%$H_{Toda}F=\varepsilon F$.
%For example, when $\g=\mathfrak{sl}_{l+1}$, the Hamiltonian
%and the eigenvalue are given by
%\bea
%&&H_{Toda}=\sum_{i=0}^lD_{i+1}D_i^{-1}(1-y_i)
%\prod_{j=i+1}^l\left(q^{-1}z_j\right),
%\quad
%\varepsilon=\sum_{i=0}^l\prod_{j=i+1}^l\left(q^{-1}z_j\right),
%\label{TodaHamiltonian}
%\ena
%where $D_i$ denotes the shift operator
%$(D_i f)(y_1,\ldots,y_l)=
%f(y_1,\ldots,qy_i,\ldots,y_l)$,
%with $D_0=D_{l+1}=1,\,y_0=y_{l+1}=0$.
%In this subsection
%we give an account on this point for completeness.
%The argument is similar to the one for
%Theorem \ref{Frec}.
%In place of the Drinfeld Casimir element, we use
%other central elements including the quadratic Casimir element.

The following construction is standard.
\begin{prop}\cite{FRT}
Let $\pi_V:\Uv\to \End(V)$
be a finite dimensional representation,  and set
$\Theta_V=(\id\otimes\pi_V)(\Theta)$,
$\Theta'_V=(\id\otimes\pi_V)(\Theta_{21})$
 and $v^{\rho_V}=\id\otimes \pi_V(v^\rho)$.
Denote further by $\varphi_V$ the operator on
$\V^\la\otimes V$ which acts
as $\id\otimes \pi_V(v^{\la-\beta})$ on each
$(\V^\la)_\beta\otimes V$.
Then, for any $k\in \Z$, the operator on $\V^\la$ given by
\bea
\mathcal{C}^{(k)}_V={\rm tr}_V\left(v^{2\rho_V}
(\Theta'_V\circ \varphi_V^{-1} \circ
\Theta_V \circ\varphi_V^{-1})^k
\right)
\label{quadCas}
\ena
acts as a scalar.
\end{prop}

In the following we take $k=-1$
and $\pi_V$ to be the vector representation for
the series $A_l,B_l,C_l,D_l$.
In terms of orthogonal vectors $\e_i$,
the simple roots are given by
\be
&&\al_1=\e_1-\e_2,\ \al_2=\e_2-\e_3,\cdots,
\al_l=\e_l-\e_{l+1}
\quad \mbox{ for $A_l$ },
\\
&&\al_1=\e_1-\e_2,\ \al_2=\e_2-\e_3,\cdots,
\al_{l-1}=\e_{l-1}-\e_{l},
\al_l=\e_l
\quad \mbox{ for $B_l$ },
\\
&&\al_1=\e_1-\e_2,\ \al_2=\e_2-\e_3,\cdots,
\al_{l-1}=\e_{l-1}-\e_{l},
\al_l=2\e_l
\quad \mbox{ for $C_l$ },
\\
&&\al_1=\e_1-\e_2,\ \al_2=\e_2-\e_3,\cdots,
\al_{l-1}=\e_{l-1}-\e_l,\ \al_l=\e_{l-1}+\e_l
\quad \mbox{ for $D_l$ }.
\en
We have $(\epsilon_i,\epsilon_i)=2$ for $B_l$ and 
$(\epsilon_i,\epsilon_i)=1$ in the other cases.
Unlike the fermionic
recursion, the Toda recursion depends on the choice of the
orientation.
Here we give the formulas for the standard
orientation compatible with the order $\al_1<\cdots<\al_l$.

\begin{prop}\label{TodaRec}
For algebras of type $A_l$ and $C_l$, we have the recursion relation
\begin{align*}
\Bigl(
{\rm tr}_V\,q^{\lambda+\rho-\beta}-{\rm tr}_V\,q^{\lambda+\rho}
\Bigr)J^\lambda_\beta
=\sum_{i=1}^lv^{-d_i}
{\rm tr}_V\Bigl(q^{\lambda+\rho-\beta}v^{\alpha_i}E_iF_i\Bigr)
J^\lambda_{\beta-\alpha_i}\,.
\end{align*}
For algebras of type $B_l$ and $D_l$, the same recursion relation
holds wherein the right hand side has an 
additional term
\begin{align*}
&-\frac{1-v^2}{1+v^2}{\rm tr}_V\Bigl(q^{\lambda+\rho-\beta}v^{2\alpha_l}E_l^2F_l^2\Bigr)J^\lambda_{\beta-2\alpha_l}
\quad \text{for $B_l$}\,,
\\
&-v^{-2}{\rm tr}_V\Bigl(q^{\lambda+\rho-\beta}v^{\alpha_{l-1}+\alpha_l}
E_{l-1}E_lF_{l-1}F_l\Bigr)J^\lambda_{\beta-\alpha_{l-1}-\alpha_l}
\quad \text{for $D_l$}\,,
\end{align*}
respectively.
\end{prop}
\begin{proof}
The scalar $\mathcal{C}^{(-1)}_V$ can be evaluated
on the highest weight vector $\one^\la$, giving
\be
(\mathcal{C}^{(-1)}_V\theta^\la_\beta,\bar\theta^\la_\beta)
={\rm tr}_V(v^{2(\la+\rho)})\
(\theta^\la_\beta,\bar\theta^\la_\beta).
\en
On the other hand, inserting \eqref{Rprod} into
$\mathcal{C}^{(-1)}_V$, we obtain
a sum of terms comprising scalar products
\be
((F^{<}_{\ga_1})^{m_1}\cdots (F^{<}_{\ga_t})^{m_t}
(E^{<}_{\ga_1})^{n_1}\cdots (E^{<}_{\ga_t})^{n_t}
\theta^\la_\beta, \
\bar\theta^\la_\beta),
\en
where $\ga_1<\cdots<\ga_t$ are the positive roots.
From the rules \eqref{root-vec2}, \eqref{root-vec3}
we observe that $\sigma(F^{<}_\gamma)$ is proportional to
$\bar E^{<}_\gamma$, and hence kills $\bar\theta^\la$ if
$\ga$ is non-simple.
Therefore, for both $\theta^\la_\beta$ and
$\bar \theta^\la_\beta$ we need to retain only root vectors
corresponding to the simple roots.
Then we must have $m_i=n_i$, so we obtain
\be
&&(\mathcal{C}^{(-1)}_V\theta^\la_\beta,\bar\theta^\la_\beta)
=\sum_{n_1,\cdots,n_i\ge 0}
\prod_{i=1}^l
%\frac{v_i^{-2n_i^2+4n_i}(1-v_i^2)^{4n_i}}
%{(v_i^2;v_i^2)^2_{n_i}}\ v^{-(\ga,\la-\beta+\ga)}
\frac{v_i^{2n_i^2-4n_i}(1-v_i^2)^{4n_i}}
{(v_i^2;v_i^2)^2_{n_i}}\ v^{-(\ga,\la-\beta+\ga)}
\\
&&\quad \times
{\rm tr}_V\left(v^{2(\la+\rho-\beta)+\ga}
E_1^{n_1}\cdots E_l^{n_l}F_1^{n_1}\cdots F_l^{n_l}\right)
\\
&&\quad \times
(E_1^{n_1}\cdots E_l^{n_l}\theta^\la_\beta,
\bar E_l^{n_l}\cdots \bar E_1^{n_1}\bar\theta^\la_\beta).
\en
Here we have set $\ga=\sum_{i=1}^ln_i\al_i$.

When $\pi_V$ is the vector representation, a simple check shows that the trace can be non-zero 
only when $\sum_{i=1}^ln_i\le1$; the only exceptions occur for 
types $B_l$ and $D_l$, where terms with $n_l=2$ and $n_{l-1}=n_l=1$ also contribute, respectively.

\end{proof}

Here are the more explicit expressions.
\be
&&
\sum_{i=1}^{l+1}
q^{(\la,\e_i)-i}(q^{-(\beta,\e_i)}-1)\ J^\la_\beta
=
\sum_{i=1}^{l}q^{(\la-\beta,\e_i)-i}J^\la_{\beta-\al_i}
\qquad \mbox{ for $A_l$ },
\en
\begin{align*}
&\sum_{i=1}^l\Bigl(q^{(\lambda,\epsilon_i)+2l-2i+1}\bigl(q^{-(\beta,\epsilon_i)}-1\bigr)
+q^{-(\lambda,\epsilon_i)-2l+2i-1}\bigl(q^{(\beta,\epsilon_i)}-1\bigr)
\Bigr)
J^\lambda_{\beta}
\\
&=\sum_{i=1}^{l-1}\Bigl(q^{(\lambda-\beta,\epsilon_i)+2l-2i+1}
+q^{-(\lambda-\beta,\epsilon_{i+1})-2l+2i+1}\Bigr)J^\lambda_{\beta-\alpha_i}
\\
&+(1+q)\bigl(q^{(\lambda-\beta,\epsilon_l)+1}+q^{-1}\bigr)J^\lambda_{\beta-\alpha_l}
-q^{(\lambda-\beta,\epsilon_l)+2}J^\lambda_{\beta-2\alpha_l}
\qquad\text{for $B_l$}\,,
\end{align*}
\be
&&
\sum_{i=1}^l
\left(q^{(\la,\e_i)+l-i+1}(q^{-(\beta,\e_i)}-1)
+
q^{-(\la,\e_{i})-l+i-1}(q^{(\beta,\e_i)}-1)\right)
\ J^\la_\beta
\\
&&\quad=
\sum_{i=1}^{l-1}\left(q^{(\la-\beta,\e_i)+l-i+1}
+q^{-(\la-\beta,\e_{i+1})-l+i}\right)J^\la_{\beta-\al_i}
+q^{(\la-\beta,\e_{l})+1}J^\la_{\beta-\al_l}
\\
&&\qquad\qquad\qquad\qquad\qquad\qquad\qquad\qquad\qquad\qquad
\qquad \mbox{ for $C_l$ },
\en
\be
&&
\sum_{i=1}^l
\left(q^{(\la,\e_i)+l-i}(q^{-(\beta,\e_i)}-1)
+
q^{-(\la,\e_{i})-l+i}(q^{(\beta,\e_i)}-1)\right)
\ J^\la_\beta
\\
&&\quad=
\sum_{i=1}^{l-1}
\left(q^{(\la-\beta,\e_i)+l-i}
+q^{-(\la-\beta,\e_{i+1})-l+i+1}\right)\ J^\la_{\beta-\al_i}
\\
&&\quad+
\left(q^{(\la-\beta,\e_{l-1})+1}
+q^{(\la-\beta,\e_{l})}\right)\ J^\la_{\beta-\al_l}
-q^{(\la-\beta,\e_{l-1})+1}\ J^\la_{\beta-\al_{l-1}-\al_l}
\\
&&\qquad\qquad\qquad\qquad\qquad\qquad\qquad\qquad\qquad\qquad
\qquad \mbox{ for $D_l$ }.
\en

%%%%%%%%%%%%%%%%%%%%%%%%%%%%%%%%%%%%%

\section{Functions $\cX^{\mu,\la}_\al$,
$\bcX^{\la,\mu}_\al$,  $X^{\la_1,\la_2}_\beta$
and recurrence relations}\label{sec:4}

In this section we give
an interpretation of the fermionic sum \eqref{pminfty}
in the context of Whittaker vectors,
and derive various relations among them and $J^\la_\beta$.

\subsection{Functions $\cX^{\mu,\la}_{\al}$,
$\bar\cX^{\la,\mu}_{\al}$}

Let $\la,\mu\in P$.
We are going to define rational functions in the variables
$z_i=q^{-(\la,\al_i)}$.
The definition is given under the assumption
$\la+\mu+2\rho\in -P_+$. This restriction is dropped afterwards.

If $\la+\mu+2\rho\in -P_+$, we have the decomposition into isotypic components
(see Corollary \ref{-P})
\bea
&&\V^{\mu}\otimes\V^{\la}=\oplus_{\al\in Q_+}
\W^{\mu+\la-\al},
\label{decompVerma1}
\\
&&\bar\V^{\mu}\otimes
\bar\V^{\la}=\oplus_{\al\in Q_+}
\bar\W^{\mu+\la-\al}\,,
\label{decompVerma2}
\ena
with each summand $\W^{\mu+\la-\al}$
(resp. $\bar\W^{\mu+\la-\al}$)
being isomorphic to a direct sum of
Verma modules  $\V^{\mu+\la-\al}$
(resp. $\bar\V^{\mu+\la-\al}$).

Choosing $\nu=(\nu_i)$ as in  \eqref{nij},
let us consider the decomposition of the vector
\bea
&&\one^\mu\otimes \theta^\la(\e,\nu,c)
=\sum_{\al\in Q_+}\eta^{\mu+\la-\al}(\e,\nu,c)
\quad \in \V^{\mu}\otimes\widehat{\V}^{\la}
\label{decompvec}
\ena
corresponding to \eqref{decompVerma1}.
Consider similarly
\be
&&\bar\one^\mu\otimes \bar\theta^\la(\e,\nu,c)
=\sum_{\al\in Q_+}\bar\eta^{\mu+\la-\al}(\e,\nu,c)
\quad \in \bar\V^{\mu}\otimes
\widehat{\bar\V}^{\la}.
\en
We define
\bea
\cX^{\mu,\la}_\al=
v^{-(\al,\al)/2+(\al,\lambda)}
\Bigl(
(\eta^{\mu+\la-\al}(\e,\nu,c))_0,
(\bar\eta^{\mu+\la-\al}(\e,\nu,c))_0\Bigr)\,.
\label{defcX}
\ena
Here $(\cdot)_0$ stands for the highest component.

In a similar manner, let us introduce
\be
&&\theta^\la(\e,\nu,c)\otimes \one^\mu
=\sum_{\al\in Q_+}\xi^{\la+\mu-\al}(\e,\nu,c)
\quad \in \widehat{\V}^{\la}\otimes\V^{\mu},
\\
&&\bar\theta^\la(\e,\nu,c)\otimes \bar\one^\mu
=\sum_{\al\in Q_+}\bar\xi^{\la+\mu-\al}(\e,\nu,c)
\quad \in
\widehat{\bar\V}^{\la}\otimes \bar\V^{\mu}.
\en
We define
\bea
\bcX^{\la,\mu}_\al=
v^{(\al,\al)/2-(\al,\lambda+2\rho)}
\Bigl(
(\xi^{\la+\mu-\al}(\e,\nu,c))_0,
(\bar\xi^{\la+\mu-\al}(\e,\nu,c))_0
\Bigr).
\label{defbcX}
\ena
Later it will be shown that \eqref{defcX},\eqref{defbcX}
are independent of the data $\e,\nu,c$
(see Proposition \ref{prop:XJJ} and Proposition \ref{lem:ident}).

\begin{prop}
For any $\mu\in P$ we have
\bea
J^{\la}_\beta
&=&\sum_{\al\in Q_+}
\cX^{\mu,\la}_\al J^{\mu+\la-\al}_{\beta-\al}
\label{cX}
\\
&=&\sum_{\al\in Q_+}
\bcX^{\la,\mu}_\al J^{\la+\mu-\al}_{\beta-\al}
q^{-(\al,\al)/2+(\al,\la+\mu+\rho)-(\mu,\beta)}.
\label{bcX}
\ena
\end{prop}
\begin{proof}
The vector
$\one^\mu\otimes\theta^{\la}(\e,\nu,c)$
is a joint eigenvector of %$E_iK_{\nu^-_i}$.
$E_iK_{\nu_i}$ with eigenvalue $c_i'=v^{(\mu,\nu_i+\al_i)}c_i$.
Therefore each isotypic component
$\eta^{\mu+\la-\al}(\e,\nu,c)$
is proportional to $\theta^{\mu+\la-\al}(\e,\nu,c')$,
where
the proportionality being determined by the highest component
$(\eta^{\mu+\la-\al}(\e,\nu,c))_0$.
Similarly,
$\bar\one^{\mu}\otimes\bar\eta^{\la}(\e,\nu,c)$
is a joint eigenvector of $\bar E_i\bar K_{\nu_i}$
with eigenvalue $c_i''=v^{-(\mu,\nu_i)}c_i$, so that
$\bar\eta^{\mu+\la-\al}(\e,\nu,c)$
is proportional to
$\bar\theta^{\mu+\la-\al}(\e,\nu,c'')$.
Formula \eqref{cX} follows from these facts.
The case of \eqref{bcX} is similar.
\end{proof}

\subsection{Function $X^{\la_1,\la_2}_{\beta}$}

Fix $\la_1,\la_2\in P$ and $\beta\in Q_+$.
Consider the Whittaker and dual Whittaker vectors
\be
&&\theta^{(1)}=\theta^{\la_1}(\e,\nu^{-},c^{(1)})
\quad \in\widehat{\V}^{\la_1},
\\
&&\theta^{(2)}=\theta^{\la_2}(-\e,-\nu,c^{(2)})
\quad \in\widehat{\V}^{\la_2},
\\
&&\bar\theta^{(1)}=\bar\theta^{\la_1}(\e,\nu,\bar c^{(1)})
\quad \in\widehat{\bar\V}^{\la_1},\\
&&\bar\theta^{(2)}=\bar\theta^{\la_2}(-\e,-\nu^{+},\bar c^{(2)})
\quad \in\widehat{\bar\V}^{\la_2},
\en
where $\nu^{\pm}_i=\nu_i\pm \al_i$, and
$c^{(i)},\bar c^{(i)}$ are chosen to satisfy
\bea
&&\bar c^{(1)}_i=v_i^{-1}c^{(1)}_i,\quad
\bar c^{(2)}_i=v_i^{3}c^{(2)}_i,
\label{cbar1}
\\
&&c^{(1)}_i/c^{(2)}_i
=-v^{(\nu_i,\la_1+\la_2-\beta-\al_i)+(\al_i,\al_i)}.
\label{cbar2}
\ena
Note that for all $\beta\in Q_+$ we have
\begin{lem}\label{lem:X}
The following conditions are satisfied for all $i$.
\bea
&&E_i\left(\theta^{(1)}\otimes\theta^{(2)}\right)_\beta=0,
\label{sing1}\\
&&\bar E_i\left(\bar\theta^{(1)}\otimes\bar\theta^{(2)}\right)_\beta=0.
\label{sing2}
\ena
\end{lem}
\begin{proof}
Noting that
\be
\Delta(E_i)=
v^{(\al_i,\nu_i-\al_i)}K_{-\nu_i+\al_i}E_i K_{\nu_i-\al_i}\otimes 1
+K_i\otimes v^{-(\al_i,\nu_i)}K_{\nu_i}E_iK_{-\nu_i},
\en
and using the defining relation \eqref{Whit}
for Whittaker vectors, we find
\be
&&E_i\left(\theta^{(1)}\otimes\theta^{(2)}\right)
=
\frac{1}{1-v_i^2}\left(v^{(\al_i,\nu_i-\al_i)}c^{(1)}_i+
v^{-(\al_i,\nu_i)}c^{(2)}_i\ K_{\nu_i}\otimes K_{\nu_i}\right)
\left(K_{-\nu_i+\al_i}\theta^{(1)}\otimes \theta^{(2)}\right).
\en
Therefore with the choice of \eqref{cbar2}
the condition \eqref{sing1} is satisfied.
Similarly \eqref{sing2} holds if
\be
\bar c^{(1)}_i/\bar c^{(2)}_i
=-v^{(\nu_i,\la_1+\la_2-\beta-\al_i)-(\al_i,\al_i)}.
\en
Because of \eqref{cbar1} and \eqref{cbar2},
it is also satisfied.
\end{proof}

We define
\bea
X^{\la_1,\la_2}_{\beta}
&=&
\Bigl((\theta^{(1)}\otimes\theta^{(2)})_\beta,
(\bar\theta^{(1)}\otimes\bar\theta^{(2)})_\beta\Bigr).
\label{defX}
\ena

\begin{prop}\label{prop:XJJ}
We have
\bea
X^{\la_1,\la_2}_{\beta}
&=&
\sum_{\al\in Q_+}J^{\la_1}_{\beta-\al}J^{\la_2}_{\al}
q^{(\al,\al)/2-(\al,\la_2+\rho)}.
\label{XJJ}
\ena
In particular, $X^{\la_1,\la_2}_{\beta}$
is independent of the data $\e,\nu$ and
$c^{(i)},\bar c^{(i)}$.
\end{prop}
\begin{proof}
This follows from the definition along with
the transformation laws
\eqref{trlaw1}, \eqref{trlaw2}.
\end{proof}

The relation \eqref{XJJ} is a counterpart of \eqref{X1}.
In order to describe the identification we introduce two vectors
$\la$ and $\bar\nu$ defined by
\be
z_i=q^{-(\la,\al_i)}, \quad \bar\nu=\sum_{i=1}^l \nu_i\omega_i.
\en

\begin{prop}
We have
\be
X^{\la_1,\la_2}_{\beta}=X^{(0,\nu)}_{C,m}(q,z),
\en
where $C$ is a Cartan matrix of ADE type and
\bea\label{identX}
\beta=\sum_{i=1}^lm_i\al_i,\quad \la_1=\la-\bar\nu,\quad  \la_2=\beta-\la-2\rho.
\ena
\end{prop}
\begin{proof}
Follows from  \eqref{XJJ} and \eqref{X1}.
\end{proof}

\begin{cor}
Let $\g$ be of ADE type. Then
we have
\bea
X^{\la_1,\la_2}_{\beta}=X^{\la_2,\la_1}_{\beta}.
\label{X12X21}
\ena
\end{cor}
\begin{proof}
The equality \eqref{X12X21} is obtained from \eqref{symmetry of X} substituting
\eqref{identX}.
\end{proof}

\begin{conj}
Relation \eqref{X12X21} holds for arbitrary $\g$.
\end{conj}

\begin{prop}
We have
\bea
X^{\la_1,\la_2}_{\beta}\bigl|_{q\to q^{-1}}
=X^{\la_2,\la_1}_{\beta}\cdot q^{(\beta,\rho)}.
\label{Xq1/q}
\ena
\end{prop}
\begin{proof}
The relation \eqref{Xq1/q} follows from
\eqref{XJJ} and \eqref{Jqinv}.
\end{proof}

The following vanishing property of
$X^{\la_1,\la_2}_\beta$ will play a key role in the sequel.
For $\mu\in P$, we denote by $L^\mu$ the
irreducible quotient of $\V^\mu$.

\begin{prop}\label{lem:vanish}
Let $\beta\in Q_+\backslash\{0\}$.
Assume that $-\la_1-\la_2-2\rho+\beta\in P_+$,
and that either $\la_1-\beta\in P_+$ or $\la_2-\beta\in P_+$.
Then we have
\bea
%X^{\la-\mu,\beta-\la-2\rho}_\beta=0
X^{\la_1,\la_2}_\beta=0
\quad \mbox{ if $-\la_1-\la_2-2\rho$ is not a weight of
$L^{-\la_1-\la_2-2\rho+\beta}$.}
\label{vanish}
\ena
\end{prop}
\begin{proof}
We apply Proposition \ref{Im} in Appendix \ref{app:2},
choosing $\la_1=\la-\mu$, $\la_2=\beta-\la-2\rho$ and
$\mu=-\la_1-\la_2-2\rho+\beta$.
Under our assumption
%$-\la_1-\la_2-2\rho$ is not a weight of
%$L^{-\la_1-\la_2-2\rho+\beta}$, hence
$v=(\theta^{(1)}\otimes\theta^{(2)})_\beta$
can be written as $\sum_{i=1}^l F_iv_i$.
Since $\bar v=(\bar\theta^{(1)}\otimes\bar\theta^{(2)})_\beta$ is
a singular vector,
we have
\be
X^{\la_1,\la_2}_\beta=\sum_{i=1}^l (F_iv_i,\bar v)
=\sum_{i=1}^l(v_i,\bar E_i\bar v)=0.
\en
\end{proof}

\begin{lem}
For any $\mu\in P$, the following recurrence relations hold.
\bea
X^{\la_1,\la_2}_{\beta}
&=&\sum_{\al\in Q_+}
\cX^{\mu,\la_1}_\al X^{\la_1+\mu-\al,\la_2}_{\beta-\al}
\label{XX1}
\\
&=&\sum_{\al\in Q_+}
\bar\cX^{\la_2,\mu}_\al X^{\la_1,\la_2+\mu-\al}_{\beta-\al}.
\label{XX2}
\ena
\end{lem}
\begin{proof}
This follows from substituting
\eqref{cX},\eqref{bcX} into \eqref{XJJ}.
\end{proof}

We can now state the relationship between
$\cX^{\mu,\la}_\beta$, $\bcX^{\la,\mu}_\beta$ and
$X^{\la_1,\la_2}_\beta$.
\begin{prop}\label{lem:ident} We have
\be
%X^{\la-\mu,\beta-2\rho-\la}_\beta
%=\cX^{\mu,\la}_\beta
%=\bcX^{\beta-2\rho-\la,\mu}_\beta.
X^{\la_1,\la_2}_\beta
&=& \cX^{\beta-\la_1-\la_2-2\rho,\la_1}_\beta \quad if \quad \la_2-\beta\in P_+  \\
&=& \bcX^{\la_2,\beta-\la_1-\la_2-2\rho}_\beta \quad if \quad   \la_1-\beta\in P_+.
\en
\end{prop}
\begin{proof}
We first assume $\la_2-\beta\in P_+$.
In the relation \eqref{XX1}, choose
$\mu=\beta-\la_1-\la_2-2\rho$ and apply
Proposition \ref{lem:vanish}. Then
the summand is non-zero only if $\al-\beta$
is a weight of $L^0$,
i.e., only if $\al=\beta$.
The first equality of Proposition follows from this.
Likewise the second follows from \eqref{XX2}.
\end{proof}

In summary, we obtain the following relations.
\begin{thm}\label{thm:XJ}
For any $\mu\in P$ we have
\bea
J^{\la}_\beta
&=&\sum_{\al\in Q_+}
X^{\la,\al-\mu-2\rho}_\al J^{\mu-\al}_{\beta-\al}
\label{id1}\\
&=&\sum_{\al\in Q_+}
X^{\al-\mu-2\rho,\la}_\al
J^{\mu-\al}_{\beta-\al}
q^{-(\al,\al)/2+(\al,\mu+\rho)-(\mu-\la,\beta)},
\label{id2}\\
X^{\la_1,\la_2}_\beta
&=&
\sum_{\al\in Q_+}
X^{\mu-\al,\la_2}_{\beta-\al}
X^{\la_1,\al-\mu-2\rho}_{\al}
\label{id3}\\
&=&
\sum_{\al\in Q_+}
X^{\la_1,\mu-\al}_{\beta-\al}
X^{\al-\mu-2\rho,\la_2}_{\al}.
\label{id4}
\ena
\end{thm}
\begin{proof}
By Proposition \ref{lem:ident},
Theorem is a restatement of the relations
\eqref{cX},\eqref{bcX}, \eqref{XX1}, \eqref{XX2}
 applied with shifted $\mu$.
\end{proof}
Identity \eqref{id1} corresponds to
\eqref{recursion J=JX},
while
identity \eqref{id3} corresponds to
\eqref{X=XX}.
We have thus shown
that these quasi-classical decompositions are exact
in the case where $C$ is a Cartan matrix of ADE type.

As an application we prove the following
\begin{thm}\label{thm:k=0-1}
For any $\beta\in Q_+$, we have
\bea
\sum_{\al\in Q_+}J^{\al-\la-2\rho}_\al
J^{\la-\al}_{\beta-\al} q^{-(\al,\al)/2+(\la+\rho,\al)}
=\delta_{\beta,0}.
\label{k=-1}
\ena
If $\g$ is of ADE type, we also have
\bea
\sum_{\al\in Q_+}J^{\al-\la-2\rho}_\al
J^{\la-\al}_{\beta-\al} =\frac{1}{(q)_\beta}.
%\qquad (\beta=\sum_{i=1}^lm_i\al_i).
\label{k=0}
\ena
\end{thm}
\begin{proof}
We first note that we expect \eqref{k=0} to hold for general $\g$.
The only reason we restrict ourselves to the ADE case is that the proof
of \eqref{k=0} uses Corollary \ref{X12X21}.

Substituting \eqref{XJJ} to \eqref{id1}, we obtain
\be
J^{\la-\mu}_\beta
=\sum_{\al,\ga}
J^{\la-\mu}_{\al-\ga}
J^{\beta-2\rho-\la}_{\ga}
J^{\la-\al}_{\beta-\al}\
q^{(\ga,\ga)/2-(\ga,\beta-\rho-\la)}.
\en
This is a linear relation among
$\{J^{\la-\mu}_\delta\}_{\delta\in Q_+}$ viewed as functions
of $\mu$.
Since this is a linearly independent set, we can compare the
coefficients of $J^{\la-\mu}_0$. In the right hand side only
the term with $\gamma=\al$ contributes, so that
\be
\delta_{\beta,0}=\sum_{\al}
J^{\beta-2\rho-\la}_{\al}
J^{\la-\al}_{\beta-\al}\
q^{(\al,\al)/2-(\al,\beta-\rho-\la)}.
\en
Renaming $\al$ to $\beta-\al$,  then changing
$\la$ to $\beta-2\rho-\lambda$, we arrive at
\eqref{k=-1}.

Similarly, after using the symmetry \eqref{X12X21} and
substituting \eqref{XJJ} to \eqref{id1}, we find
\be
J^\la_\beta=\sum_{\al,\ga}J^{\al-\mu-2\rho}_{\al-\ga}J^{\la}_{\ga}
q^{(\ga,\ga)/2-(\ga,\la+\rho)}J^{\mu-\al}_{\beta-\al}.
\en
Specializing $q^{-\la}$ to $0$,
the left hand side simplifies due to \eqref{Jz=0}, while
only the term with $\ga=0$ remains in the
right hand side.
This proves \eqref{k=0}.
\end{proof}

Finally we give the counterpart of the identity
\eqref{relation 2} in Section \ref{sec:2}.

\begin{thm}\label{thm:J[0,k]}
For  a non-negative integer $k$ we have
\begin{align}
\sum_{\alpha\in Q_+}
J_\alpha^{\al-\lambda-2\rho}
J^{\lambda-\alpha}_{\beta-\alpha}
q^{k((\alpha,\alpha)/2-(\lambda+\rho,\alpha))}
\label{formula 1}
=J^\la_\beta[0,k],
\end{align}
where the right hand side is defined in \eqref{eq:frec2}.
\end{thm}
\begin{proof}
Let us denote
the left hand side of \eqref{formula 1} by
$\bar J^\la_\beta[0,k]$, and set
\be
\bar J^\la_\beta[r,s]
=q^{r\left(\frac{(\beta,\beta)}2-(\lambda+\rho,\beta)\right)}
\bar J^\la_\beta[0,s-r].
\en
The previous formula \eqref{k=0} states that
$\bar J^\la_\beta[0,0]=J^\la_\beta[0,0]$.
Using \eqref{Jrs3} for $J^\la_\beta=J^\la_\beta[0,\infty)$, it is easy
to verify that
\be
\bar J^\la_\beta[0,k]=\sum_{\beta_1+\beta_2=\beta}
\bar J^{\la-\beta_2}_{\beta_1}[0,0]\bar J^\la_{\beta_2}[1,k].
\en
The same relation holds for $J^\la_\beta[0,k]$ by \eqref{Jrs3}.
Hence by induction we obtain
$\bar J^\la_\beta[0,k]=J^\la_\beta[0,k]$.
\end{proof}
We remark that \eqref{formula 1} in the limit $k\rightarrow\infty$
reproduces the fermionic formula for $J^\lambda_\beta$.

%%%%%%%%%%%%%%%%%%%%%%%%%%%%%%%%%%%%%%%%%%
\appendix

\section{Direct proof of Toda recursion for $\mathfrak{sl}_{l+1}$}
\label{app:1}

{\rm
We give here a direct proof that
the fermionic sum $I_{C,m}(q,z)$ for
the Cartan matrix $C$ of type $A_l$  satisfies the
Toda recursion.
Since we fix $C$, we drop it and
denote $I_{C,m}(q,z),W_{C,m}$ by $I_m(q,z),W_m$.
}

\begin{prop}
The rational functions $I_m(q,z)$ $(m\in\Z_{\geq0}^l)$
are characterized by the Toda recursion
\begin{align}
&{\textstyle\left\{\sum_{i=0}^l\left(q^{m_{i+1}-m_i}-1\right)\prod_{j=i+1}^l\left(q^{-1}z_j\right)\right\}}I_m(q,z)
\label{Toda}\\
&\quad=\sum_{i=1}^l{\textstyle\left\{q^{m_{i+1}-m_i}\prod_{j=i+1}^l\left(q^{-1}z_j\right)\right\}
I_{m_1,\cdots,m_i-1,\cdots,m_l}(q,z)}.\nonumber
\end{align}
Moreover, they satisfy the symmetry relation
\begin{align}
I_{m_l,\ldots,m_1}(q^{-1},z_l^{-1},\ldots,z_1^{-1})=
I_{m_1,\ldots,m_l}(q,z_1,\ldots,z_l)(qz)^mq^{W_m}.\label{symmetry}
\end{align}
\end{prop}
\begin{proof}
Let $\tilde I_m(q,z)$, $m\in \Z_{\ge 0}^l$ be a set of rational functions in
$q,z=(z_1,\ldots,z_l)$ such that $\tilde I_0(q,z)=~1$.
It is straightforward to show by induction on $m$ that for $\tilde I_m(q,z)$
the Toda recursion \eqref{Toda}
implies the symmetry relation \eqref{symmetry}.
Now, we want to show that the former also implies the
fermionic recursion \eqref{fermionic recursion}.
Since the solution is unique for both \eqref{Toda} and
\eqref{fermionic recursion},
the statement of Proposition follows.

Let $\C(q,z)$ be the field of rational functions in
$q,z=(z_1,\ldots,z_l)$.
Consider the vector space over $\C(q,z)$ consisting
of formal power series in $y=(y_1,\ldots,y_l)$
with coefficients in $\C(q,z)$. We denote it by $\mathcal F$.

We consider the $\C(q,z)$-linear actions $y_i,D_i$
$(i=1,\ldots,l)$ on $\mathcal F$:
\begin{align*}
y_i\cdot f(y_1,\ldots,y_l)&=y_if(y_1,\ldots,y_l),\\
D_i\cdot f(y_1,\ldots,y_l)&=f(y_1,\ldots,qy_i,\ldots,y_l),
\end{align*}
and set formally $D_0=D_{l+1}=1,\,y_0=y_{l+1}=0$.

Let $\tilde I_m(q,z)\in\C(q,z)$. We assume that
$\tilde I_0(q,z)=1$.
Set
$$F(q,z,y)=\sum_my^m\tilde I_m(q,z)$$
and
$$G(q,z,y)=\sum_my^m\tilde I_m(q,z)z^mq^{W_m}.$$
They belong to $\mathcal F$.
Set
$$H=\sum_{i=0}^l\left(D_{i+1}D_i^{-1}(1-y_i)-1\right)\prod_{j=i+1}^l\left(q^{-1}z_j\right).$$
The Toda recursion reads as $HF=0$, and the symmetry relation reads as
\begin{align}
F\left(q^{-1},z_l^{-1},\ldots,z_1^{-1},q^{-1}y_l,\ldots,q^{-1}y_1\right)=G(q,z,y).\label{s in F G}
\end{align}
Set
\begin{align}
\Lambda=\prod_{i=1}^l\frac1{(y_i)_\infty}=\sum_m\frac{y^m}{(q)_m}.
\end{align}
The fermionic recursion reads as
\begin{align}
F(q,z,y)=\Lambda G(q,z,y).\label{f r in F G}
\end{align}
Our goal is to show that if $HF=0$, and therefore \eqref{s in F G} is valid, then \eqref{f r in F G} follows.

Suppose that $HF=0$. By changing $q\rightarrow q^{-1},z_i\rightarrow z_i^{-1},y_i\rightarrow q^{-1}y_i,
D_i\rightarrow D_i^{-1}$, we obtain
$$
\sum_{i=0}^l\left(D_{i+1}^{-1}D_i(1-q^{-1}y_i)-1\right)\prod_{j=i+1}^l\left(qz_j^{-1}\right)
F\left(q^{-1},z_1^{-1},\ldots,z_l^{-1},q^{-1}y_1,\ldots,q^{-1}y_l\right)=0.
$$
Using
$$\Lambda^{-1}D_i\Lambda=(1-y_i)D_i,$$
we can rewrite this as
\begin{align}
\left(\Lambda^{-1}H\Lambda\right)F\left(q^{-1},z_l^{-1},\ldots,z_1^{-1},q^{-1}y_l,\ldots,q^{-1}y_1\right)=0.
\label{Lambda}
\end{align}
Because of the uniqueness of the solution $HF(q,z,y)=0$ with $F(q,z,0)=1$,
we obtain
$$F(q,z,y)=\Lambda F\left(q^{-1},z_l^{-1},\ldots,z_1^{-1},q^{-1}y_l,\ldots,q^{-1}y_1\right).$$
From \eqref{s in F G} and this equality follows the fermionic recursion \eqref{f r in F G}.
\end{proof}

%%%%%%%%%%%%%%%%%%%%%%%%%%%%%%%%%%%%%%%%%%%

\section{Proposition on singular vectors}\label{app:2}

{\rm
The main goal of this Appendix is to prove
a statement about
singular vectors which is used in the main text.
In what follows,
for a $\Uv$ module $M$,
$[M]_\nu$ will denote its subspace of weight $\nu$.
}

We start with the following Lemma.

\begin{lem}\label{sum}
Let $M$ be a $\Uv$ module from the category $\mathcal{O}$.
Let $p\in [M]_{-\mu-2\rho}$, $\mu\in P_+$
be a singular vector such that
\be
p\notin \sum_{i=1}^l \mathrm{Im} F_i.
\en
Then
the Verma module $\Uv\cdot p$ generated by $p$ is a direct summand
in $M$.
\end{lem}
\begin{proof}
We first note that since the Verma module $\V^{-\mu-2\rho}$ is irreducible,
the submodule $V=\Uv\cdot p$ generated by $p$ is isomorphic to it.
We now show that there exists
a submodule $W\subset M$ such that
$M=V\oplus W.$

Denoting by $C_v$ the quantum Drinfeld Casimir element
$v^{\sum_{i=1}^l\al_i^{\vee}\omega_i+2\rho}u$
(see \cite{Dri} and Proposition \ref{DCas}) we have the
decomposition
$M=\oplus_{z\in \K} M^z$
into the generalized eigenspaces
\be
M^z=\{ m\in M\mid (C_v-z)^km=0\quad \mbox{ for some $k$ }\}.
\en
Setting $z_0=v^{(\mu,\mu+2\rho)}$
we have $V\subset M^{z_0}$ and
$M=M^{z_0}\oplus \bigoplus_{z\neq z_0}M^{z}$.
So it suffices to find a submodule $M_0\hk M^{z_0}$
such that $M^{z_0}=V\oplus M_0$.

Since $M\in\mathcal{O}$, there exists a sequence of submodules
\be
L_0\hk L_1\hk L_2\hk\dots,\quad L_j\hk M^{z_0}/V,\quad \lim_{j\to\infty} L_j=M^{z_0}/V
\en
such that each quotient $L_j/L_{j-1}$ is a highest weight module with highest weight $\tau_j$.
For all $j\ge 0$,
let $\bar w_j\in L_j$ be a vector of weight $\tau_j$ such that
the image of $\bar w_j$ in $L_j/L_{j-1}$ is highest weight
vector.
Note  that $\tau_j\not< -\mu-2\rho$. In fact,
for all $j$ we have
$z_0=v^{(\tau_j,\tau_j+2\rho)}$.
If we set $\al= -\mu-2\rho-\tau_j$, we obtain
\be
(\al,\al)+2(\rho,\al)+2(\mu,\al)=0,
\en
which
implies that $\tau_j\not\in -\mu-2\rho-(Q_+\backslash\{0\})$.
Hence, if $\tau_j\ne -\mu-2\rho$, there exists the unique vector $w_j\in M^{z_0}$
which is a lifting of $\bar w_j$.
For $j$, such that $\tau_j= -\mu-2\rho$, we fix arbitrary liftings
$w_j\in M^{z_0}$. Setting
\be
M_0= \sum_{j\ge 0} \Uv\cdot w_j,
\en
we obtain $V + M_0= M^{z_0}$.
Since
$V=\Uv\cdot p$ is irreducible and $p\notin \sum_{i=1}^l \mathrm{Im} F_i$,
the intersection $V \cap M_0$ is trivial.
This proves the Lemma.
\end{proof}

\begin{cor}\label{-P}
Let $M$ be a $\Uv$ module from the category $\mathcal{O}$ such that
$M=\bigoplus_{\mu\in P_+} [M]_{-\mu-2\rho}$.
Then $M$ is isomorphic to a direct sum of Verma modules.
\end{cor}
\begin{proof}
Let $M_0\subset M$ be the maximal submodule such that $M=M_0+W$ is a decomposition
into the direct sum of $\Uv$ modules and $M_0$ is a direct sum of Verma modules.
Let $w\in W$ be a singular vector such that $[W]_\la=0$ for $\la$ bigger
than the weight of $w$. Then Lemma $\ref{sum}$ implies $W=(\Uv\cdot w)\oplus W'$
and thus $M_0$ is not maximal.
\end{proof}

\begin{prop}\label{Im}
Let $L^\mu$ be an irreducible $\Uv$ module
with highest weight $\mu\in P_+$,
and let $\la\in P$, $\beta\in Q_+$
be such that either
$\la+2\rho\in -P_+ $ or
$\mu+\beta-\la\in -P_+$.
Assume further that $[L^\mu]_{\mu-\beta}=0$.
Then we have
\be
\left[\V^{\la-\mu}\T \V^{\beta-2\rho-\la}\right]^{sing}_{-\mu-2\rho}
\subset
\sum_{i=1}^l \mathrm{Im} F_i
\en
where $(~)^{sing}$ means the space of singular vectors.
\end{prop}

\begin{proof}
Set $M=\V^{\la-\mu}\T \V^{\beta-2\rho-\la}$, and
suppose that the statement of the
Proposition is not true.
Then there exists a vector
$p\in M^{sing}_{-\mu-2\rho}$
such that $p\notin \sum_{i=1}^l \mathrm{Im} F_i$.
Set $V=\Uv\cdot p$.
Because of Lemma \eqref{sum}, there exists
a submodule $W\subset M$ such that
\bea\label{W}
M=V\oplus W.
\ena
Tensoring both sides of \eqref{W} by $L^\mu$
we obtain
\bea\label{VL}
\V^{\la-\mu}\T \V^{\beta-2\rho-\la}\T L^\mu=
(\V^{-\mu-2\rho}\T L^\mu) \oplus (W\T L^\mu).
\ena
We show that the decomposition \eqref{VL} is impossible
by a homological argument.

In the following we set $U=\Uv$.
Let $\na$ (resp. $\bo$, $H$) be the subalgebra of $U$
generated by $\{E_i\}_{1\le i\le l}$ (resp.
$\{E_i,K_i^{\pm 1}\}_{1\le i\le l}$,
$\{K_i^{\pm 1}\}_{1\le i\le l}$). All these subalgebras are
vector spaces over the field $\K=\C(v^{1/\mathcal{N}})$.
We shall make use of the following facts.

($i$)
Let $X$ be a $\bo$ module and
$
\mathrm{Ind}_{\bo}^U X=U\T_{\bo} X
$
be the induced $U$ module.
Then we have an isomorphism
\bea\label{TorInd}
\Tor_\bullet^{U,H}(\K,\mathrm{Ind}_{\bo}^U X)
\simeq \Tor_\bullet^{\bo,H}(\K, X).
\ena
The proof is essentially given in \cite{K}, Lemma $3.1.14$, which treats the classical
case of \eqref{TorInd}. In order to prove \eqref{TorInd} we only need to
replace the classical
$(\bo, H)$ projective resolution of $X$ from \cite{K}, Corollary $3.1.8$ by an
arbitrary $\bo$-free resolution.

($ii$) Denoting by $n$ the number of positive roots of $\g$
we have
\bea
&&
\dim [\Tor^\na_n(\K,\V^\la)]_{\nu}
=\delta_{\nu,\la+2\rho},
\label{Tor1}
\\
&&\dim
[\Tor^{\na}_n(\K,L^\mu)]_{\nu}=\delta_{\nu,\mu+2\rho}.
\label{Tor2}
\ena
The proof of these equalities is based on the quantum BGG resolution
\cite{HK} (see also \cite{M}, \cite{R}), which generalizes the classical
BGG resolution  \cite{BGG}.
Let $W$ be the Weyl group of $\g$. For  $w\in W$ and $\la\in P$ we denote by
$l(w)$ the length  of $w$ and by $w*\la=w(\la+\rho)-\rho$ the shifted action of $w$.
In order to prove \eqref{Tor1} we use
the quantum BGG resolution of the trivial $\Uvv$ module
\bea\label{BGGK}
0\to F_n\to F_{n-1}\to \cdots\to F_0\to\K\to 0,
\ena
where $F_p=\bigoplus_{w\in W: l(w)=p} \bV^{w*0}$ is a direct sum of Verma modules $\bV^{w*0}$
over $\Uvv$. Using the anti-isomorphism $\sigma:\Uv\to \Uvv$ \eqref{sigma} we endow
each $F_p$ with the structure of right $\Uv$ module. Thus \eqref{BGGK} is a right
$\bo$-free resolution of the trivial module $\K$ and
$\Tor^\na_n(\K,\V^\la)$ is equal to $n$-th homology of the complex
\bea\label{nth}
0\to F_n\T_{N} \V^\la\to F_{n-1}\T_{N} \V^\la\to\ \cdots\to F_0\T_{N} \V^\la\to 0.
\ena
We note that $F_n$ is the free $\bo$ module with one generator of $H$-weight
$-w_0*0$, where $w_0$ is the longest element in $W$. Since $\V^\la$ is irreducible
we obtain that the space of $n$-th homology of \eqref{nth} is one-dimensional and is
generated by the tensor product of highest weight vectors of $\bV^{w_0*0}$ and of $\V^\la$.
Now the equality $w_0\rho=-\rho$ implies \eqref{Tor1}. The proof of \eqref{Tor2}
is very similar and uses the (left) quantum BGG resolution of the module $L^\mu$.

For any $\la$ we have
$\V^{\la}=\mathrm{Ind}^U_{\bo} \K_\la$,
where $\K_\la$ denotes the one-dimensional $\bo$ module
with trivial action of $\na$
and an action of $K_i$ by $v^{(\la,\al_i)}$.
This gives (see \cite{K}, Proposition $3.1.10$)
\be
\V^{\la-\mu}\T \V^{\beta-2\rho-\la}\T L^\mu
&=&
\mathrm{Ind}_{\bo}^U (\K_{\la-\mu}\T
\V^{\beta-2\rho-\la}\T L^\mu) \\
&=&
\mathrm{Ind}_{\bo}^U (\K_{\beta-2\rho-\la}\T
\V^{\la-\mu}\T L^\mu).
\en
We conclude that
\bea
\Tor_n^{U,H}(\K,\V^{\la-\mu}\T \V^{\beta-2\rho-\la}\T L^\mu)
&=& \left[
\Tor_n^\na(\K,\V^{\beta-2\rho-\la}\T L^\mu)\right]_{\mu-\la}
\label{Tor3a} \\
&=& \left[
\Tor_n^\na(\K,\V^{\la-\mu}\T L^\mu)\right]_{\la+2\rho-\beta}
\label{Tor3b}.
\ena
Suppose $\mu+\beta-\la\in -P_+$.
Then we have the decomposition (see Corollary \ref{-P})
\be
\V^{\beta-2\rho-\la}\T L^\mu=\oplus_{\nu}
\V^{\beta-\la-2\rho+\nu}\otimes [L^\mu]_{\nu}.
\en
From \eqref{Tor1} and the vanishing assumption,
the right hand side of \eqref{Tor3a} is equal to
$\left[L^\mu\right]_{\mu-\beta}=0$.
Now suppose
$\la+2\rho\in -P_+$.
Then
we have the decomposition (see Corollary \ref{-P})
\be
\V^{\la-\mu}\T L^\mu=\oplus_{\nu}
\V^{\la-\mu+\nu}\otimes [L^\mu]_{\nu}.
\en
Again, from \eqref{Tor1} and the vanishing assumption,
the right hand side of \eqref{Tor3b} is equal to
$\left[L^\mu\right]_{\mu-\beta}=0.$

Similarly, \eqref{Tor2} implies
\be
\Tor_n^{U,H}(\K,\V^{-\mu-2\rho}\T L^\mu)
\simeq \left[\Tor_n^{\na}(\K, L^\mu)\right]_{2\rho+\mu}
=\K.
\en
This shows that the decomposition \eqref{W} is impossible,
and thus proves our Proposition.
\end{proof}

%%%%%%%%%%%%%%%%%%%%%%%%%%%%%%%%%%%%%%%%%%%

\bigskip

{\it Acknowledgement}.\quad
Research of BF is partially supported by RFBR grants 08-01-00720-a,
NSh-3472.2008.2 and 07-01-92214-CNRSL-a.
Research of EF is supported by the
RFBR  Grants
06-01-00037, 07-02-00799 and NSh-3472.2008.2, by Pierre Deligne fund based on his 2004
Balzan prize  in mathematics, by Euler foundation and by Alexander von Humboldt
Fellowship.
Research of MJ is supported by the Grant-in-Aid for Scientific
Research B-20340027 and B-20340011.
Research of TM is supported by
the Grant-in-Aid for Scientific Research B--17340038.
Research of EM is supported by NSF grant DMS-0601005.
The present work has been carried out during the visits
of BF, EF and EM to Kyoto University.
They wish to thank the University for hospitality.
This work was supported by World Premier International
Research Center Initiative (WPI Initiative), MEXT, Japan.

The authors thank A. Tsymbaliuk for pointing out that a term was missing
in the Toda Hamiltonian of type $B_l$.

\end{document}